\documentclass[preprint,12pt]{article}

\usepackage{amsmath,amsfonts,amsthm,amssymb,graphicx,color}
 \usepackage{nonfloat}
 \usepackage{amsmath}
 \usepackage{amsfonts}
 \usepackage{amssymb}
 \usepackage{graphicx}
\usepackage{url}
\usepackage[dvips]{epsfig}		
\usepackage{array}						
\usepackage{multicol}					
\usepackage{xspace}						
\usepackage{indentfirst} 			
\usepackage{latexsym}					
\usepackage{setspace}					
\usepackage{makeidx}					
\usepackage{vmargin}					
\usepackage{enumerate}
\usepackage{ucs}
 \usepackage{times}
\usepackage{epsfig}
\usepackage{epstopdf}
\usepackage[active]{srcltx}
\usepackage[mathscr]{eucal}
\usepackage{caption}

\usepackage{ucs}
\usepackage{epsfig}
\usepackage{epstopdf}
\usepackage[active]{srcltx}
\usepackage[mathscr]{eucal}
\usepackage{caption}
\usepackage{subfigure}
\usepackage{tabularx}
\usepackage{caption}
\usepackage{float}
\restylefloat{figure}
\usepackage{subfig}
\usepackage{tikz}
\usepackage{xcolor}
\def\b{\beta}
\def\a{\alpha}
\def\g{\gamma}
             
\def\as{A_s}
\def\bi{B_{int}}
\def\ba{B_{abs}}
\makeindex

\newcommand{\Bo}{\textbf}
\newcommand{\insoluble}{%
\draw (0,0) -- (4,0);
\draw [ color=black ] ( 0,-.7 ) ellipse (0.2cm and 0.7cm)  ;
\draw [ color=gray!40, fill=gray!40  , thin ] ( 1,-.5 ) circle (0.2cm) ;
\draw [ color=gray!40, fill=gray!40  , thin ] ( .4,-.3 ) circle (0.2cm) ;
\draw [ color=gray!40, fill=gray!40  , thin ] ( .4,-1.1 ) circle (0.2cm) ;
\draw [ color=gray!40, fill=gray!40  , thin ] ( 3.4,-.4 ) circle (0.2cm) ;
\draw [ color=gray!40, fill=gray!40  , thin ] ( 0.6,-1 ) circle (0.2cm) ;
\draw [ color=gray!40, fill=gray!40  , thin ] ( 1.6,-0.5 ) circle (0.2cm) ;
\draw [ color=gray!40, fill=gray!40  , thin ] ( 2,-1 ) circle (0.2cm) ;
\draw [ color=gray!40, fill=gray!40  , thin ] ( 2.3,-1 ) circle (0.2cm) ;
\draw [ color=gray!40, fill=gray!40  , thin ] ( 2.1,-0.8 ) circle (0.2cm) ;
\draw [ color=gray!40, fill=gray!40  , thin ] ( 2.4,-.3 ) circle (0.2cm) ;
\draw [ color=gray!40, fill=gray!40  , thin ] (3.4,-1.1 ) circle (0.2cm) ;
\draw [ color=gray!40, fill=gray!40  , thin ] (3.8,-.9 ) circle (0.2cm) ;
\draw [ color=gray!40, fill=gray!40  , thin ] (1.5,-1.0 ) circle (0.2cm) ;
\draw [ color=black ] ( 4,-.7 ) ellipse (0.2cm and 0.7cm)  ;
\draw (0,-1.4) -- (4,-1.4);
\draw (0,-.7 ) -- (-0.2,-.9) node[left]{\small $r$};
}
\newtheorem{definition}{Definition}[section]

\def\b{\beta}
\def\a{\alpha}
\def\g{\gamma}
             
\def\as{A_s}
\def\bi{B_{int}}
\def\ba{B_{abs}}

\begin{document}
\title{Digestion Modelling in the Small Intestine : Impact of Dietary Fibre}

\author{M. Taghipoor$^{\dagger \star}$, G. Barles\footnote{
Laboratoire de Math\'ematiques et Physique Th\'eorique (UMR CNRS 6083).
F\'ed\'eration Denis Poisson (FR CNRS 2964) Universit\'e de Tours.
Facult\'e des Sciences et Techniques, Parc de Grandmont, 37200
Tours, France}, C. Georgelin$^{\star}$,\\ J.R. Licois$^{\star}$
\& Ph. Lescoat\footnote{INRA, UR83 Recherches Avicoles, 37380 Nouzilly, France.}}
\date{}
\maketitle

\begin{abstract}
In this work, we continue the modelling of the digestion in the small intestine, started in a previous article, by investigating the effects of dietary fibre. We recall that this model aims at taking into account the three main phenomena of the digestion, namely the transit of the bolus, the degradation of feedstuffs and the absorption through the intestinal wall. 
In order to study the role of dietary fibre on digestion, we model their two principal physiochemical characteristics which interact with the function of the small intestine, i.e. viscosity and water holding capacity. This leads us to consider some features of digestion which have not been taken into account previously, in particular the interrelationship between the evolution of dry matter and water in the bolus. The numerical results  are in agreement with the positive effect of insoluble dietary fibre on the velocity of bolus along the small intestine and on its degradation. These results highlight the negative effect of soluble dietary fibre on digestion. Although, this model is generic and contains a large number of parameters, to our knowledge, it is among the first qualitative dynamical modelling of fibre influence on intestinal digestion.
\end{abstract}

\vspace{.5cm} {\small \noindent{\bf Key-words :}  Modelling of intestinal digestion, water holding capacity, soluble dietary fibre, insoluble dietary fibre.

\section{Introduction}\label{Intro}

Digestion in the small intestine can be described through three main phenomena~: transit of the bolus along the small intestine, degradation of macromolecules into smaller ones and absorption through intestinal wall. Taking into account these phenomena, the authors have presented in \cite{Taghipoor2011} a generic model of digestion in which the bolus include only one category of macromolecules (carbohydrates, proteins or lipids) and water.

However mixing these nutrients influences the digestion process through interactions between molecules.
In order to improve this model and to make it more realistic, we should consider the effects of such interactions, and we have decided to do so by first including dietary fibre in the bolus because of their significant role on the digestion. One of the key properties of fibre is its water holding capacity and this leads us to investigate the role of water in the digestion processes. To do so, we distinguish dry matter and water in each substrate, we model water kinetic in correlation with the dry matter one and we take into account the impact of water in all the aspects of digestion.

To be more precise on our approach, we first describe the main underlying assumptions which guide our modelling, then we show how they are translated into equations and finally numerical tests are performed for examining the effects of several hypothesis. In addition, we point out known or assumed mechanisms relating dietary fibre and digestion process.

This paper is organized as follows : physiochemical characteristics of soluble and insoluble dietary fibre are introduced in Section \ref{WDF}. In Section~\ref{key}, main assumptions of the model are presented. Section~\ref{prel} is devoted to the description of the composition of bolus, the chemical transformations of macromolecules and all the notations. Digestion in the presence of dietary fibre includes the modification of transport as well as degradation and absorption equations as described in Section~\ref{Model}. Section~\ref{num} is a comparison of the numerical results considered as \textit{in silico} experiments. Finally, in Section~\ref{disc}, the model is discussed and perspectives are proposed.\\

\noindent{{\bf Acknowledgement.}
The multidisciplinary collaboration on this research project between the INRA Center of Nouzilly and the Laboratoire de Math\'ematiques et Physique Th\'eorique was initiated within and supported by the CaSciModOT program (CAlcul SCIentifique et MOD\'elisation des universit\'es d'Orl\'eans~et de Tours) which is now a Cluster of the french Region Centre.
This collaboration also takes place in a CNRS-INRA PEPS program ``Compr\'ehension et Mod\'elisation du devenir de l'aliment dans le tube digestif``. This work is part of the PhD thesis of Masoomeh Taghipoor, financed by CNRS and INRA.}

\section{Biological Background on Water and Dietary Fibre}\label{WDF}
Digestion modelling requires the knowledges of the physiochemical properties of macromolecules concerned by this phenomenon as well as mechanical and biochemical reactions observed for their degradation.

Dietary fibre (DF) is usually defined as the sum of plant non-starch polysaccharides and lignin that are not hydrolysed by the enzymes secreted by the non-ruminant digestive system, but that can be partially digested by microflora in the gut. A main effect of fibre is to regulate intestinal degradation and absorption of nutrients as well as their transit along the gut. Physiochemical characteristics of fibre include viscosity, hydration, fermentability (mostly in the large intestine), adsorption or entrapment of nutrients and bulking effect. Each of these characteristics affects meaningfully the function of the  gastrointestinal tract \cite{Wenk2001,McCleary2001}. These characteristics depend on the polysaccharides chemistry. One way to classify dietary fibre is based on their water solubility. Insoluble dietary fibre include cellulose, some hemicelluloses and lignin. The other is soluble dietary fibre such as viscous fibre which includes beta-glucans, pectins, gums, mucilages and some hemicelluloses \cite{Anderson2009,Tharakan2009}.

For monogastrics, most available nutrients are degraded and absorbed in the small intestine. At the beginning of duodenum the bolus consists of partially degraded feedstuffs and water.  Once in the small intestine, mechanical and chemical digestion of feedstuffs make the nutrients available to the organism. Enzymatic hydrolysis is the most important chemical reaction in digestion, which takes place in aqueous solution.  Enough water is required for an efficient digestion process eventhough water/nutrient ratio are not precisely known. Furthermore, classification of dietary fibre through their water solubility and the impact of Water Holding Capacity (WHC) of DF on digestion reveals the key-role of water on digestion. WHC is defined by B. Shneeman \cite{McCleary2001} as the ability of fibre source to swell when mixed with water and to hold water within its matrix. 
 
 \subsubsection*{Soluble dietary fibre}\label{SF}
Soluble DF are believed to impact significantly digestion and absorption as well as bolus transport in the small intestine. The main physiochemical properties of soluble DF are viscosity, water holding capacity (WHC) and organic compound entrapment \cite{K.E.2001}.  Soluble DF, because of its high viscosity, is generally associated with slow transit through the stomach and increasing of the small intestinal transit time \cite{Siljestroem1986}.

\subsubsection*{Insoluble dietary fibre}
Insoluble DF acts primarily in the large intestine where, due to its WHC, increases faecal bulk, dilutes colonic contents and decreases mouth-to-anus transit time \cite{K.E.2001}. However, its effects on digestion and transit in the small intestine can not be neglected since insoluble DF affects the transit time in the small intestine through its laxative property. Recent studies have shown that the inclusion of a moderate level of dietary fibre improves the digestibility in chicks \cite{Jim'enez-Moreno2009}. Therefore to obtain an optimal efficiency in nutrient utilization, Burhalter et al.  \cite{Burkhalter2001} proposed to increase the ratio of insoluble to soluble DF. Moreover, the use of insoluble fibre in commercial broiler chicks improves the intestine morphological parameters and result in a better performance assumed to be connected to more efficient digestion and absorption processes \cite{Rezaei2011}.
Two hypothesis are proposed in order to study the influence of insoluble DF on nutrients digestibility in the small intestine : ($i$) insoluble DF increases the retention time in the stomach changing the nutrient profile of the bolus at the entry of the small intestine which could lead to a higher digestion and absorption. ($ii$) Physical characteristics of insoluble DF change the digestion process mostly through their capacity of swelling water and nutrients. These both hypothesis are either tested in the\textit{ in silico} experiments (cf. Section~\ref{testf}) or included in the equations (cf. Section \ref{Model}).

\subsection*{Water}
To have a better understanding of the role of dietary fibre, it is therefore required to study more precisely the evolution of water during digestion in the small intestine. For example, depending on the bolus composition, water absorption through intestinal wall or its secretion into the small intestine lumen could be observed. However, the evolution of water amount in the small intestine depends also on other components' kinetics within the bolus as described by (\ref{degAs}), (\ref{degbi}) and (\ref{degba}) in Section \ref{prel}.
%
%
%
\section{ Key Model Assumptions}\label{key}
In this section, the key assumptions for the model are presented. 
\begin{itemize}
 \item[\textbf{H1:}]\textit{Each component of the bolus (macromolecules, partially degraded macromolecules, nutrients and fibre) is represented mathematically as a portion of dry matter and a characteristic proportion of water.}
\end{itemize}
For example, ``starch in a bolus'' includes both dry starch and water used to maintain starch molecules in aqueous solution. The same is observed for the ``disaccharides in a bolus'' and ``glucose in a bolus'' combining smaller molecules resulting from starch hydrolysis associated with a specific level of water. In other words, a component $C$ of bolus is represented as $C^{dm}+W_C$ where $C^{dm}$ denotes the dry matter of $C$ and $W_C$ is the necessary amount of water to maintain it in a solution state. Moreover, we have assumed that the mass of $W_C$ is proportional to the mass of $C^{dm}$, i.e. equal to $c\, C^{dm}$ for some characteristic number $c\geq0$ whic represents the necessary amount of water proportion to maintain $C$ in solubilized phase. Despite the presence of water in the bolus, a little amount of non-solubilized dry matter may be included in bolus, which is (of course) associated with $c=0$. 
\begin{itemize}
 \item[\textbf{H2:}] \textit{Without DF, the bolus contains a single macromolecule and water. It is represented by a homogeneous cylinder with the constant length $\ell$. Including insoluble DF transforms this homogeneous bolus into an heterogeneous one by modifying the concentrations of feedstuffs and nutrients.}
\end{itemize}
A bolus in the small intestine is a viscous solution of dry matter and water: we assume that 
 its volume is very close to the volume occupied by the water in the bolus. In other words, our assumption is that $DM$ does not fill any volume (e.g. : solubilized sugar + water does fill the same volume as the water alone).    
Once the volume is known, [H2] allows to compute the radius of the bolus $r(t)$, which impacts its movement along the small intestine. For a bolus with constant mass at the beginning of duodenum, [H2] influences the degradation by increasing the concentration of nutrients in the homogeneous part of bolus. Indeed the space occupied by insoluble dietary fibre is unavailable for the other macromolecules.
\begin{itemize}
 \item[\textbf{H3:}]\textit{Digestion in the small intestine is due to volumic and surfacic transformations. Volumic degradation is the enzymatic hydrolysis of bolus components by pancreatic and exogenous enzymes inside the bolus while the surfacic one is the degradation by brush border enzymes on the internal wall of small intestine.}
\end{itemize}
Some additional facts should be pointed out : water facilitates the contact of the macromolecules with the brush border enzymes, enhancing the surfacic degradation. Increasing the water to dry matter ratio dilutes the bolus and decreases the volumic degradation. Both of these reactions are proportional to the mass ratio of concerned substrates.
\begin{itemize}
 \item[\textbf{H4:}]\textit{The bolus movement along the small intestine is due to peristaltic waves. The efficiency of these waves are proportional to the radius of bolus, and inversely proportional to the distance from the pylorus.}
\end{itemize}
 We use, in this article, the model of transport with an averaged velocity function as presented in \cite{Taghipoor2011b}. The movement of the bolus in the model by this equation is due to a homogenized acceleration caused by the average effect of the peristaltic waves. Of course, there would no difficulty to come back to the model which takes into account all the frequent peristaltic waves.

\begin{itemize}
 \item[\textbf{H5:}] \textit{The water in the bolus which is not hold by the macromolecules and DF through WHC i.e. ``available water'', decreases the viscosity of the bolus and facilitates its movement. Due to osmotic type equilibrium, the concentration of this ``available water'' tends to reach a fixed ratio.} 
\end{itemize}
 In other words, ``Available water'' reduces the friction caused by the bolus contact with the intestinal wall.
\begin{itemize}
 \item[\textbf{H6:}] \textit{Dietary fibre modifies the bolus evolution through its WHC by holding the water in its matrix, and therefore modifying the volume of the bolus. Soluble DF decreases the efficiency of peristaltic waves.}
\end{itemize}
 By their WHC, dietary fibre holds a significant quantity of water in the bolus and therefore keep the volume and the radius of the bolus higher. It is worth pointing out  that soluble dietary fibre change the consistence of the bolus in the sense of making it more jelly, implying the decreases of efficiency of peristaltic waves.

\section{Physiological Aspects and Bolus Composition}\label{prel}
Different steps of mechanical and chemical transformations are detailed below.
We present also the composition of the bolus and their interactions.

\begin{itemize}
 \item\Bo{Non degradable substrate $A_{nd} $ :} The quantity $A_{nd}$ represents the mass of macromolecules which is not degradable by endogenous enzymes of the digestive tract.

 \item\Bo{Non solubilized substrate $A_{ns}$ :} The quantity $A_{ns}$ represents the mass of macromolecules which is not accessible to enzymatic hydrolysis. In presence of a sufficient quantity of water, $A_{ns}$ is transformed into $A_s$.

 \item\Bo{Solubilized substrate $A^{dm}_s $ :} The quantity $A^{dm}_s$ is the mass of dry substrate in solution state. It is called solubilized substrate and it is assumed that one unit of $A^{dm}_s$ requires $W_s$ units of water to remain solubilized. Recalling [H1],  $W_s$ represents the required mass of water to solubilize $\as^{dm}$.  This quantity depends on the properties of each macromolecules. The mass of $W_{s}$ is assumed to be equal to the mass of $\a \as^{dm}$ where $\a$ represents the ratio of water associated with $A_s^{dm}$. For simplification purposes, the mix of $A^{dm}_s$ and water is represented  by $A_s$ .

 \item\Bo{Intermediate substrate $B^{dm}_{int}$ :} The quantity $B^{dm}_{int}$ is the mass of dry intermediate substrate obtained from the degradation of $A_s$ by volumic transformation [H3]. It is solubilized and $W_{int}$ represents the required amount of water to maintain solubilization. The mass of $W_{int}$ is assumed to be equal to the mass of $\b \bi^{dm}$ where $\b$ represents the ratio of water associated with $\bi^{dm}$. For $\bi=B^{dm}_{int}+W_{int}$, volumic transformation is represented as  
				\begin{equation}\label{degAs}
				A_s+enzymes\rightarrow B_{int} + (W_{s}-W_{int})\tag{Reaction 1}
				\end{equation}
 Depending on the value of $W_{s}$ and $W_{int}$, the amount of $W_{s}-W_{int}$ of water can be released or hold in the bolus. 

 \item\Bo{Absorbable nutrients $\ba^{dm} $ :} The quantity $\ba^{dm}$ is the mass of dry absorbable nutrients obtained from surfacic reactions (cf. [H3]). For $\ba=\ba^{dm}+ W_{abs}$, the surfacic transformation is defined as   
				\begin{equation}\label{degbi}
				\as +enzymes\rightarrow \ba+(W_{s}-W_{abs}) \tag{Reaction 2}
				\end{equation}
				\begin{equation}\label{degba}
				\bi+enzymes\rightarrow \ba + (W_{int}-W_{abs})\tag{Reaction 3}
				\end{equation}
$W_{abs}$ is the required amount of water to maintain solubilization and its mass is assumed to be equal to $\g\,\ba^{dm}$ where $\g$ represents the ratio of water associated with the $\ba^{dm}$.  

 \item\Bo{Soluble and insoluble dietary fibre :}$ F^{dm}_{sol}$ and $F^{dm}_{insol}$ represent the soluble and non-soluble dry dietary fibre respectively. The main property of dietary fibre presented in Section \ref{Intro} and hypothesis [H5] is its water holding capacity.
$$F_{sol}=W_{sol}+F^{dm}_{sol} $$
$$F_{insol}=W_{insol}+F^{dm}_{insol} $$
where the mass of $W_{sol}\,(W_{insol})$ is assumed to be equal to $\lambda_sF^{dm}_{sol} \,(\lambda_{i}F^{dm}_{insol})$ for $\lambda_s$ and $\lambda_{i}$ which represent the ratio of water associated with $ F^{dm}_{sol} $ and $F^{dm}_{insol}$ respectively . As described in Section \ref{Intro}, DF is not degradable by endogenous enzymes of the small intestine. 
 
The following diagram shows the different transformations inside the bolus.
\begin{figure}[h]
 \centering%
\fbox{\includegraphics[width=110mm]{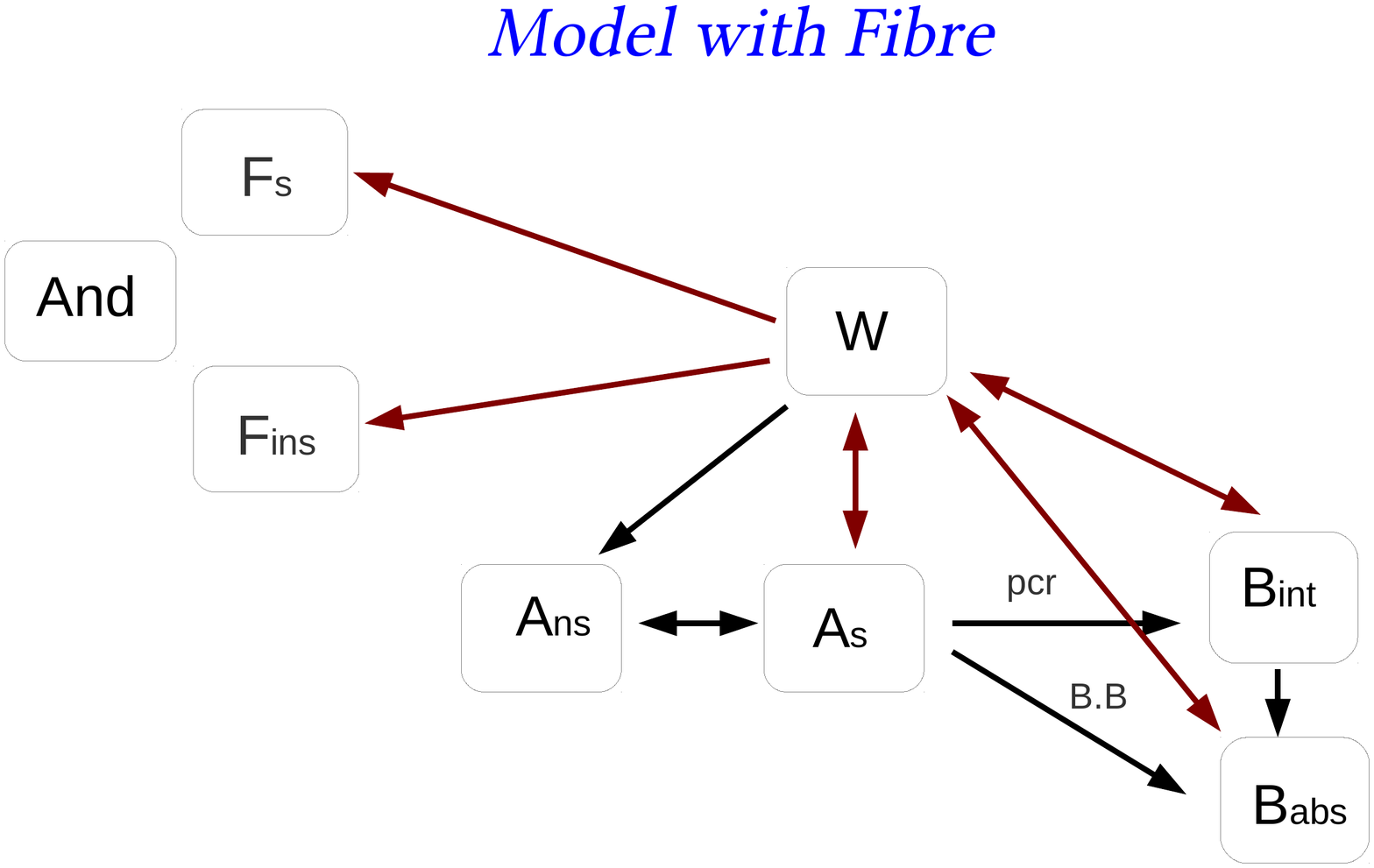}}
\caption{\small Physical and chemical transformations inside the bolus and included in the model are represented in this scheme. ``B.B`` stands for brush border enzymes and ``pcr'' stands for pancreatic ones.}\label{diagram}
\end{figure}

 \item\Bo{Dry Matter :} Total amount of Dry Matter substrate in the bolus is therefore
$$DM=A_{nd}+A_{ns} +\as^{dm}+ \bi^{dm}+\ba^{dm} + F^{dm}_{sol} +F^{dm}_{insol}.$$ 

 \item\Bo{Water :} Impact of dietary fibre on digestion is closely linked to their WHC capacity. Though, water evolution in the bolus has to be described to understand effects of DF on digestion.  

 Total water $W_{tot}$ in the bolus comes from three main sources
\begin{itemize}
 \item[ $(i)$] $W_{feed}$ : water incorporated naturally in feedstuffs (e.g. : one gram of wheat contains 12 $\%$ of water). The amount of $W_{feed}$ coming out of the stomach is assumed to be proportional to the ingested dry matter $DM$,
 $W_{feed}=K_{feed}\cdot DM$.
\item[ $(ii)$]   $W_{sec}$ : Water included in the endogenous secretions of saliva and stomach which is also assumed to be proportional to the ingested dry matter, $W_{sec}=K_{sec}\cdot DM$. 
\item[$(iii)$] 
 $W_{drink}$ : Drunk water is assumed to be independent to the quantity of $DM$. 
\end{itemize}
Total water included in the bolus is therefore defined as the sum of $W_{feed}$, $W_{sec}$ and $W_{drink}$
$$W_{tot}=W_{feed}+W_{sec}+W_{drink}.$$
Thereby, ``available water'', $W$, as presented in [H6], is defined as the difference between $W_{tot}$ and the water associated with $DM$ for maintaining the solution state. In term of mass, the quantity of ``available water'' in the bolus at each time is 
$$W= W_{tot}-W_{s}-W_{int}-W_{abs} -W_{sol}-W_{insol}\; .$$

 \item\Bo{Mass of bolus $M$ :} The total mass of bolus, $M$, is given by $M=DM+W_{tot}$.

 \item\Bo{Volume of bolus $V$ :} To define the bolus volume, as explained in [H2], the volume of each substrate in solution is assumed to be the same as the volume filled by water associated with that substrate, i.e. the volume of $A_s$ is equal to $W_s/\rho_w$, where $\rho_w=1$ is the water density.  
The volume of bolus is therefore represented as
			      \begin{equation*}\label{mass}
                              V=W_{tot}=\pi r^2\ell\rho_w=\pi r^2\ell.
                                \end{equation*}
Since the length of the bolus is assumed to be fixed, the volume evolution leads to compute the radius $r(t)$ of bolus at each time. Consequently its surface is written as $S=2\pi r\ell$.

In the following section, the different properties of dietary fibre on the model of digestion are taken into account.
\end{itemize}

\section{Model Equations}\label{Model}
To include WHC property of soluble DF in the digestion model, the mix of $F^{dm}_{sol}$ and $W_{sol}$ is assumed to form a viscous gel in the bolus. Therefore in our model the mass of ``available water'' in the bolus is reduced. Moreover, viscous fibre enhances motility but decreases transit rate, since it resists propulsive contractions \cite{FAO1998}. This resistance to peristaltic waves is described through a new notion called efficient radius of bolus called $r_{sol}$. As described in [H2], the volume filled by soluble DF is $W_{sol}=\lambda_s F_{sol}$, then 

\begin{definition}\label{Rsol}
The efficient radius of bolus is defined as
 $$ r_{sol}=\sqrt{(W_{tot}-W_{sol})/2\pi \ell}.$$
In the same way, the efficient surface of bolus is described as $S_{sol}~=~2\pi r_{sol}\ell$.
\end{definition}
This definition is used to describe the decrease of the surfacic degradation and absorption caused by soluble DF in the model.      

One of the hypothesis in the first model of digestion in \cite{Taghipoor2011} is the bolus homogeneity. The mass concentration of each component of the bolus is assumed to be its mass divided by the total mass $M$ of the bolus.
To model the digestion in presence of insoluble DF, new notions are defined because of heterogeneity of bolus as described in [H2].
The volume filled by insoluble DF (the mix of $F_{insol}^{dm}$ and $W_{insol} $), is assumed to be unavailable to the macromolecules of feedstuffs in the bolus. In the digestion model, this hypothesis is taken into account by the following definition. 

\begin{definition}\label{conc} The apparent concentration of different substrates in the bolus is represented as  
$$[\as^{dm}]=\dfrac{\as^{dm}}{M-F_{insol}  }\,,\quad [\bi^{dm}]=\dfrac{\bi^{dm}}{M-F_{insol}  }\,,[\ba^{dm}]=\dfrac{\ba^{dm}}{M-F_{insol}  }$$
$$[A_{ns}]=\dfrac{A_{ns}}{M-F_{insol} }\,,\quad  [W]=\dfrac{W}{M-F_{insol} }\,,\quad [A_{nd}]=\dfrac{A_{nd}}{M-F_{insol}  }$$
\end{definition}

The degradation of macromolecules $\as$ and $\bi$ as well as the absorption of the nutrients $\ba$ are affected by WHC property of the insoluble DF through this definition. 
Figure \ref{insoluble} shows the regions of bolus which are filled by insoluble DF and therefore unreachable by the macromolecules and nutrients. 
	\vspace{1cm}
	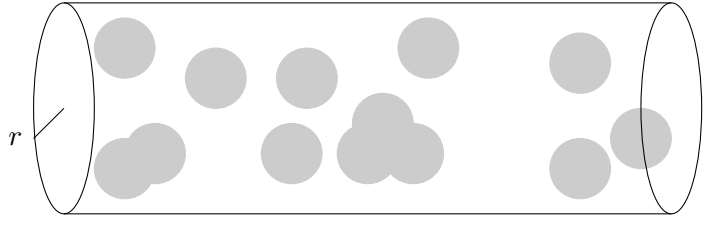
\begin{figure}[h]
	    \centering
	    \begin{tikzpicture}[scale=2]
		\insoluble
	    \end{tikzpicture}
	\caption{\small The distribution of insoluble fibre in the bolus as assumed in the model. The apparent volume $V_{app}$ (see Definition~\ref{vapp}) is the white part of the cylinder.}\label{insoluble}
	\end{figure}
	\vspace{1cm}

Integrating the insoluble DF in the bolus changes also the region reachable to volumic degradation (see Figure \ref{insoluble}). According to hypothesis [H2], $W_{insol}$ is the volume filled by insoluble fibre. 

\begin{definition}\label{vapp}
 The apparent volume of bolus $V_{app}$ is therefore defined as
$$V_{app} = W_{tot}-W_{insol}.$$ 
\end{definition}

The above considerations are taken into account in the following steps of digestion described below. 

\subsection{Transport of bolus}\label{transport}
The averaged equation of transport of bolus introduced in \cite{Taghipoor2011b} reads
	\begin{eqnarray*}
	\dfrac{d^2x}{dt^2}=\tau(1-c^{-1}\frac{dx}{dt})\frac{c_0+c_1r}{a+bx}-\dfrac{K_{visco}}{[W]}\frac{dx}{dt},
	\end{eqnarray*}
where $\tau $ is the mean effect of the pulses by unit of time, $c$ is mean velocity of peristaltic waves, $x$ the position in the small intestine. Taking into account the properties of dietary fibre, this equation changes to
	\begin{eqnarray}\label{TransSolDF}
	\dfrac{d^2x}{dt^2}=\tau(1-c^{-1}\frac{dx}{dt})\frac{c_0+c_1r_{sol}}{a+bx}-\dfrac{K_{visco}}{[W]}\frac{dx}{dt}.
	\end{eqnarray}
The bolus movement described by this equation depends on its position in the small intestine and on its efficient radius. Moreover, the acceleration is slowed down by a viscosity term which depends on the available water.
\subsection{Volumic transformation}\label{vol}
Volumic transformation presented by (\ref{degAs}) in Section \ref{prel} is the degradation of $\as^{dm}$ due to pancreatic and exogenous enzymes resulting in production of $\bi^{dm}$.
 Evolution of $\as^{dm}$ by this transformation is represented by 
	\begin{equation*}\label{as-vol}
	\dfrac{d\as^{dm}}{dt}=-k_{vol}(x)[\as^{dm}]V_{app}
	\end{equation*}
where $k_{vol}(x)$ takes into account the enzymatic activity which is a function of bolus position at each time $t$. The term $[\as^{dm}]V_{app}$ describes the dependence of volumic degradation on the concentration of $\as^{dm}$ at each unit of apparent volume i.e. the volume filled by the insoluble DF is not accessible to the enzymes and macromolecules.  

Consequently, integrating the insoluble fibre in the bolus influences the volumic transformation by increasing the substrates concentration via the Definition \ref{conc} and by changing the volume and using $V_{app}$ introduced in Definition \ref{vapp}.

As described in Section \ref{transport}, integrating the soluble fibre modifies the velocity of bolus along the small intestine and therefore the distance travelled at each time $x(t)$. Consequently, it influences the volumic degradation through the function $k_{vol}(x(t))$. 
     
This reaction takes place in a solution and each unit of $\as^{dm}$ degraded by this equation causes the release of ``available water'' $W_s=\a\as^{dm}$. 
 
The volumic production of intermediate substrate  $\bi^{dm}$ is the result of degradation of $\as^{dm}$
	\begin{equation*}\label{bi-vol}
	\dfrac{d{B}^{dm}_{int}}{dt}=k_{vol}(x)[\as^{dm}]V_{app}.
	\end{equation*}
Each unit of produced $\bi^{dm}$ requires and uses $W_{int}=\b\bi^{dm}$ to maintain the solution state.
According to (\ref{degAs}), the result of volumic transformation is the consumption or release of ``available water'' $W$. Thereby, volumic evolution of water is represented as
	\begin{equation}\label{volW}
	\dfrac{d{W}}{dt}=\tilde{k}_{vol}(x)[\as^{dm}]V_{app} \; ,
	\end{equation} 
where $\tilde{k}_{vol}=(\a-\b) k_{vol}.$

Soluble fibre can be hydrolysed by exogenous enzymes $e_{exo}$ (cellulase, hemicellulases, ...) ingested by food
		\begin{equation}\label{Fsreac}
		F_{sol}+e_{exo} \rightarrow \bi +(W_{ol}-W_{int}).\tag{Reaction 4}
		\end{equation}

Hydrolysis of $F_{sol}$ by exogenous enzymes follows the same evolution as the volumic transformation of $\as^{dm}$
		\begin{equation*}\label{Fs}
		\dfrac{d{F}_{sol}^{dm}}{dt}= -k_{s}e_{exo}\tilde{ph}(x)[F_{sol}^{dm}]V_{app}
		\end{equation*}
where $[F^{dm}_{sol}]=F_{sol}/(M-F_{insol})$ and $\tilde{ph}(x)$ is the exogenous enzyme activity along the small intestine.  
This reaction produces the intermediate substrate $\bi^{dm}$ 

	\begin{equation*}\label{bi-vol2}
	\dfrac{d{B}^{dm}_{int}}{dt}=k_{s}e_{exo}\tilde{ph}(x)[F_{sol}^{dm}]V_{app}.
	\end{equation*}
The amount $(W_{sol}-W_{int})$ is released by (\ref{Fsreac}) and modifies the evolution of water 
	\begin{equation}\label{FsW}
	\dfrac{d{W}}{dt}=...+\tilde{k}_s e_{exo}\tilde{ph}(x)[F_{sol}^{dm}]V_{app}
	\end{equation} 
where $\tilde{k}_s=(\lambda_s-\beta)k_{s}$.
\subsection{Surfacic transformation}\label{surf}
Surfacic degradation is the last step of transformation of macromolecules in the small intestine. The produced nutrients by this degradation are then absorbed through intestinal wall. 
Surfacic degradation depends on the fraction of $\as$ on the surface of the bolus represented by  $\displaystyle[\as]S_{sol}$ therefore 
	\begin{equation*}\label{As-sur}
	\dfrac{d{A}^{dm}_{s}}{dt}= -k_{surf}[\as][W]S_{sol}
	\end{equation*}
where $k_{surf}$ stands for the rate of surfacic degradation of $\as$ and the efficient surface $S_{sol}$ has been defined by Definition \ref{Rsol}. Moreover, it is assumed that the brush-border enzymes are always in excess in the small intestine.
Surfacic degradation of $\bi^{dm}$ follows the same process as for $\as^{dm}$. Therefore for $[\bi]$ defined by Definition \ref{conc}, we have 
	\begin{equation*}\label{bi-sur}
	\dfrac{d{B}^{dm}_{int}}{dt}= -\tilde{k}_{surf}[\bi][W]S_{sol}.
	\end{equation*}
where $k_{surf}$ stands for the rate of surfacic degradation of $\bi$.
Evolution of water in the bolus is influenced by the surfacic degradation, i.e. the quantity of water consumed (or released) by (\ref{degbi}). Therefore the surfacic evolution of water is
	\begin{equation}\label{surfW}
	\dfrac{d{W}}{dt}=...+\bigg((\b-\g)\tilde{k}_{surf}[\bi]+ (\a-\g)k_{surf}[\as]\bigg)[W]S_{sol}.\end{equation} 
\subsection{The equilibrium between $\as $ and $A_{ns} $ }\label{equil}
Modifications of feedstuffs in the stomach by the enzymes and water change most of the $A_{ns}$ into $A_s=A_s^{dm}+W_s$ and makes them accessible to intestinal enzymes. 

However, for some feedstuffs, the bolus may contain $A_{ns}$ at the beginning of small intestine. In this case, the digestion of bolus contains also the transformation of $A_{ns}$ into $A_s$. Mixing with bile acid for lipids and producing the micelles, denaturing for the proteins and adding water and solubilization for the dry starch are examples of the transformation of $A_{ns}$ into $\as$ in the small intestine.

The solubilization of $A_{ns}$ which takes place in the presence of enough quantity of $W$ and results in the production of $\as$, is a phenomenon taken into account in the model. Solubilization is reversible and lack of water may cause production of $A_{ns}$ releasing $W$ in the bolus. 

Thereby, the balance is assumed to be reached when
 $$A_s=\mu([W])A_{ns}$$
for $\mu$ which is an increasing function of $[W]$. If $k_{equi}$ stands for the rate of turning back to equilibrium then the dynamical equilibrium may be defined 
	\begin{equation}\label{equilAns}
\dfrac{d{A}_{ns}}{dt}=-k_{equi}\big(\mu([W])A_{ns}-\as\big)
	\end{equation}
and therefore 
	\begin{equation*}\label{equilAs}
	\dfrac{d{A}_{s}^{dm}}{dt}=k_{equi}\big(\mu([W])A_{ns}-\as\big).
	\end{equation*}
The variation of water quantity caused by the equilibrium may be represented as 
	\begin{eqnarray}\label{equiW}
	&\dfrac{d{W}}{dt}=...+\a k_{equi}\big(\mu([W])A_{ns}-\as\big).
	\end{eqnarray}
\subsection{Pancreatic and biliary secretions}\label{panc}
Pancreatic and biliary secretions consist of a solution of nutrients and enzymes which do not contain available water $W$. In fact, water included in this solution is assumed to be associated with nutrients and enzymes to keep them solubilized. Modelling details on these secretions could be seen in \cite{Taghipoor2011}.
Adding dietary fibre increases the quantity of pancreatic secretions. However this point is not yet included in the model.  
\subsection{Absorption through intestinal wall}\label{abs}
Absorption of nutrients through intestinal wall depends on their concentration on the inner surface of intestinal wall 
	\begin{equation*}\label{equabs}
	\dfrac{d{B}^{dm}_{abs}}{dt}=...-k_{abs}[\ba^{dm}]S_{sol}
	\end{equation*}
The passage of nutrient through intestinal wall releases the associated water, thus
\begin{eqnarray}\label{equW}
	&\dfrac{d{W}}{dt}=...+\g k_{abs}[\ba^{dm}]S_{sol}
	\end{eqnarray}
where $k_{abs}$ represents the rate of absorption through intestinal wall.
\subsection{Water equilibrium}\label{W}

Water equilibrium was already taken into account in \cite{Taghipoor2011}. The assumption was that $[W]$ tends to reach a fixed ratio ($10\%$), suggesting the equation 
				\begin{equation}\label{watere}
				\dfrac{d[W]}{dt}=-k_w([W]-0.1) 
				\end{equation}
where $[W]=W(t)/M(t)$, $M(t)$ representing the bolus mass.
The superposition of the equations (\ref{volW}), (\ref{FsW}), (\ref{surfW}), (\ref{equiW}), (\ref{equW}) and (\ref{watere}) provides the equation describing the evolution of $W$ along the small intestine. 

The evolution of bolus mass is represented  as 
	\begin{equation*}
	\dfrac{d{M}}{dt}= (\a+1)\dfrac{d{A}_{s}^{dm}}{dt}+ (\b+1)\dfrac{d{B}_{int}^{dm}}{dt}+(\g+1)\dfrac{d{B}_{abs}^{dm}}{dt}+(\lambda_{sol}+1) \dfrac{dF^{dm}_{sol}}{dt}+\dfrac{d{A}_{ns}}{dt}+\dfrac{d{W}}{dt}, 
	 \end{equation*}
each term of the above equation is replaced by its expression, therefore we obtain
	\begin{equation}\label{equmass}
	\dfrac{d{M}}{dt}= \dfrac{M}{M-W}(-k_w(W-0.1M)-  k_{abs}[\ba^{dm}]S_{sol}).
	\end{equation}
The variation of bolus volume depends on the absorption or secretion of ``available water'' and endogenous secretions in the small intestine i.e. 
	\begin{equation}\label{equvol}
	\dfrac{d{V}}{dt}= \dfrac{d{W}}{dt},
	\end{equation}
therefore 
	    \begin{equation*}
	      \dfrac{d{ V}}{dt}= -k_w(W-0.1M)+ \hbox{Secretions}.
	    \end{equation*}
%
%
%
 \section{Numerical Simulations}\label{num}
 A thorough examination of the effects of the different parameters of the model on transport, degradation and absorption is carried out by Scilab software. In these \textit{in silico} experiences, the presence of dietary fibre $F_{sol}$ and $F_{insol}$ in the bolus, the variation of the initial values of bolus and the sensitivity of the model on the parameters $\a$, $\b$ and $\g$ (cf. Section \ref{prel}) are investigated.

The mass of the bolus at the entry of the small intestine is assumed to be fixed in all our following experiences.
The bolus at the entry of the small intestine contains $A_{nd}$, $A_{ns}$, $\as$, $\bi$, $\ba$ and $W$. When studying the influence of dietary fibre on digestion, the non degradable substrate $A_{nd}$ is replaced by $F_{sol}$ or $F_{insol}$.
%
\subsection{Influence of dietary fibre on intestinal absorption}\label{testf}
As described in Section \ref{Intro}, the positive effect of insoluble DF in digestion may be due to two main reasons : $(i)$ modification of the composition of bolus due to the increase of retention time in the stomach, $(ii)$ modification of the bolus physical characteristics.

%
The effects of these two cases on the digestion model are studied in this section.
\subsubsection*{Influence of the modification of the bolus in the stomach}
Including insoluble DF in the bolus delays gastric emptying. The direct effect of this phenomenon is to increase the solubilization in the bolus and to start partially the degradation. We studied the effects of such a change in the initial conditions for our model.  
To this aim, two numerical experiences are carried out : $(a)$ the increase in the ratio of $\as$ to $A_{ns}$ and $(b)$ the increase in the ratio of $\bi$ to $\as$ when $A_{ns}=0$. 

\begin{itemize}
 \item[$(a)$]Our first experience consists in increasing the ratio of $\as$ to $A_{ns}$ in the bolus at the entry of the small intestine. The value of absorbed dry nutrients at the end of the small intestine does not vary meaningfully. Table \ref{AsAnstab} shows the absorbed dry nutrients to $DM$ ratio at the end of ileum $x=17$ for different for the different $\as^{dm}$ to $DM$ ratio at the beginning of the small intestine.
	\begin{table}[h]
	\centering%
	\begin{tabular}{|c|c|c|}
	  \hline
	\multicolumn{2}{|c|}{$x=0$} &\multicolumn{1}{c|}{ End of the small intestine $x=17$}  \\
	\hline
	$\as^{dm}.DM^{-1}$ (\%) & $A_{ns}.DM^{-1}$ (\%) & Absorbed dry nutrients to $DM$ ratio (\%)\\
	\hline
	0& 85 &56 \\
	42&42& 57 \\
	85&0 & 58\\
	\hline
	\end{tabular}\caption{\small The proportion of absorbed nutrients to $DM$ for different scenarios of solubilization in the bolus at the entry of duodenum.}\label{AsAnstab}
	\end{table}
Despite the variations in the ratio of $\as$ to $A_{ns}$, the equilibrium between $\as$ and $A_{ns}$ defined by Equation (\ref{equilAns}) is reached quickly (see Figure \ref{AsAnsfig}). The sensitivity analysis shows also that the value of absorbed dry nutrients is almost independent to the variations in the ratio of $\as$ to $A_{ns}$. However, this results depend on the choice of equilibrium rate $k_{equi}$, a small value of $k_{equi}$ may decrease the difference between the final values of absorbed nutrients. 

	\begin{figure}[h]
	\subfigure[$\as^{dm}=36$, $A_{ns}=0$]{
	\includegraphics[width=70mm]{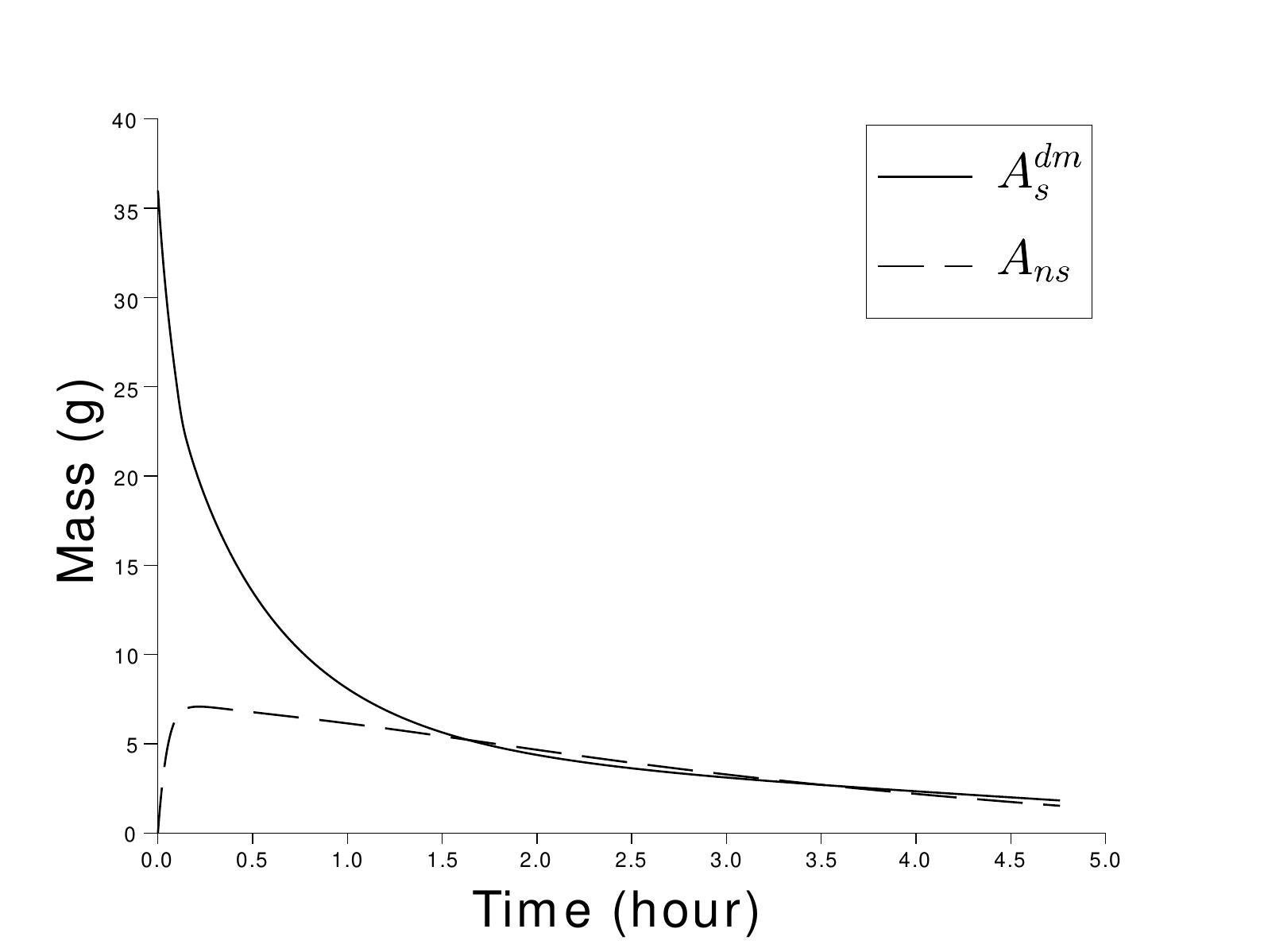} 
	  }
	\subfigure[$\as^{dm}=18$, $A_{ns}=18$]{
	\includegraphics[width=70mm]{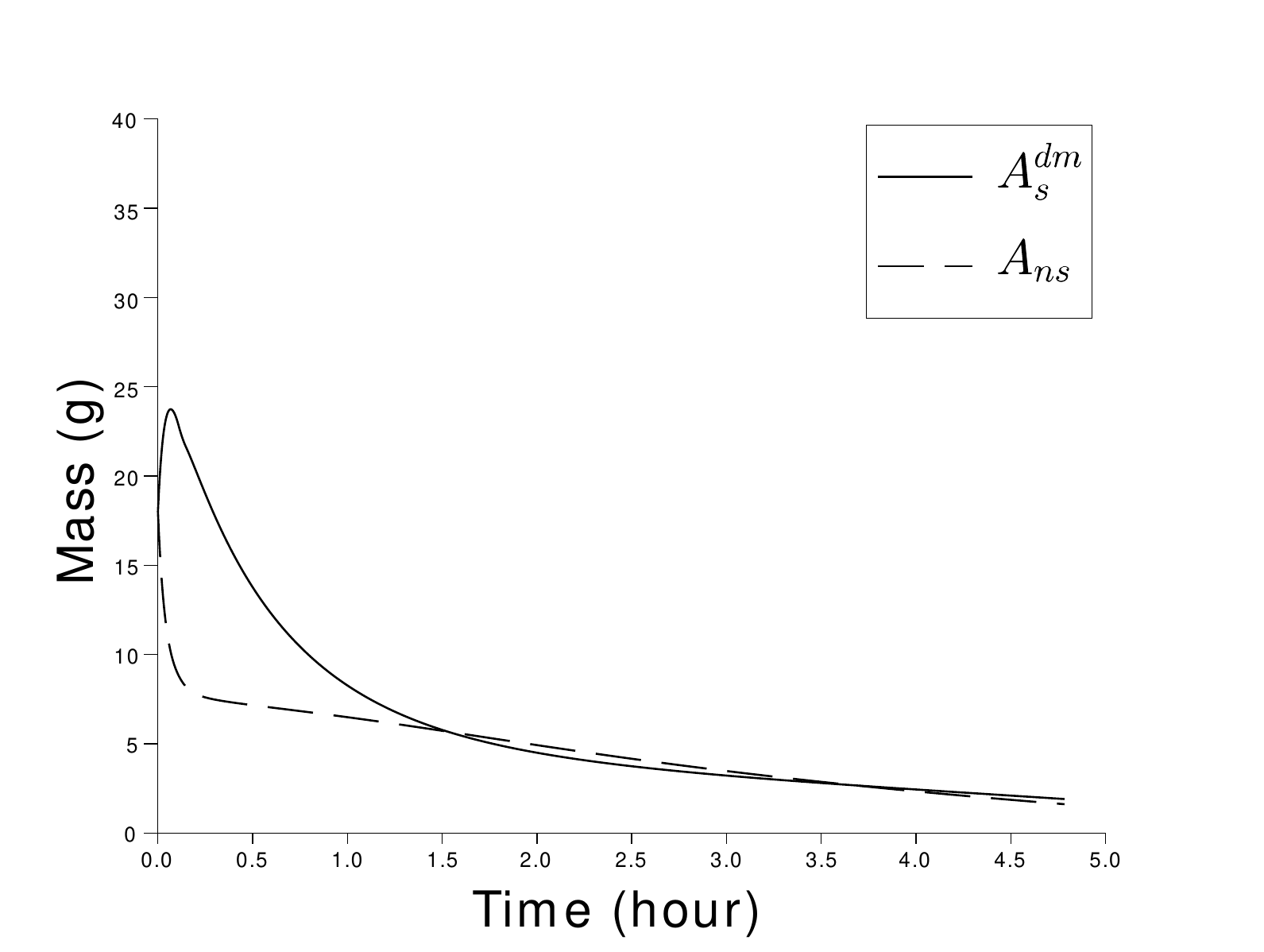} 
	  }
	\subfigure[$\as^{dm}=0$, $A_{ns}=36$]{
	\includegraphics[width=70mm]{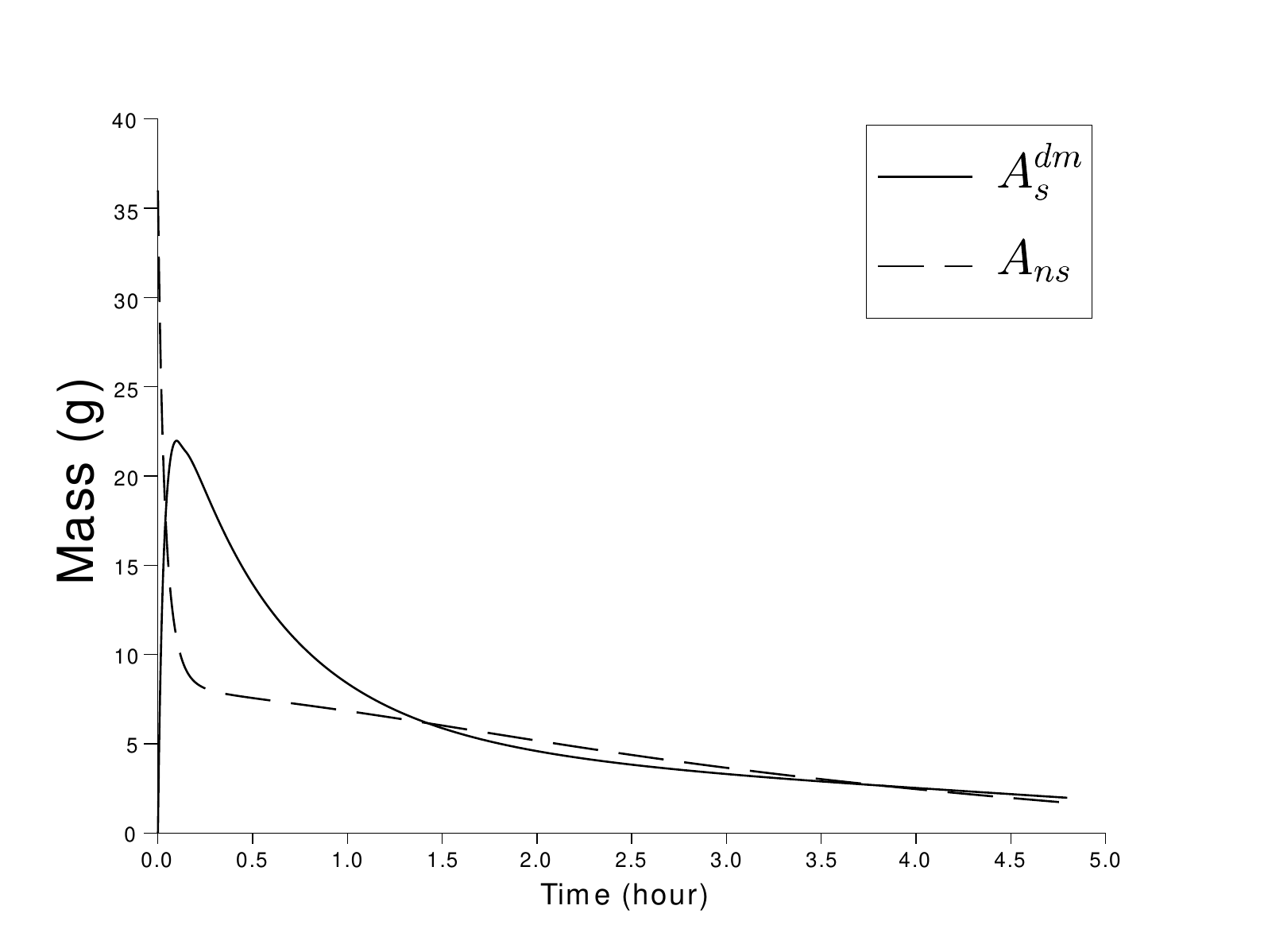}
	  }
	\caption{\small The equilibrium $\as^{dm}-A_{ns}$ is reached quickly for different initial value of $\as$.}\label{AsAnsfig}
	\end{figure}
\item[$(b)$]
In the second experience, the modification in the stomach is assumed to result in the transformation of all $A_{ns}$ in the bolus into $\as$ and additionally the production of $\bi$. Variations in the ratio of $\as$ to $B_{int}$ inside the bolus at the entry of the small intestine are tested. Table \ref{AsBint} shows the variation of absorbed dry nutrients at the end of the small intestine as a function of the initial value of $B_{int}$.

Numerical results shows the increase in absorbed dry nutrients when the ratio of $\bi$ to $\as$ is increased.
	\begin{table}[h]
	\centering%
	\begin{tabular}{|c|c|c|c|}
	  \hline
	\multicolumn{2}{|c|}{$x=0$} &\multicolumn{1}{c|}{ End of the small intestine $x=17$}  \\
	\hline
	$\bi^{dm}.DM^{-1}$ (\%) & $A_{s}.DM^{-1}$ (\%) & Absorbed dry nutrients to $DM$ ratio (\%)\\
	\hline
	0& 85 &57 \\
	42&42& 61 \\
	85&0 & 64\\
	\hline
	\end{tabular}\caption{\small The relation between the absorbed dry nutrients at the end of digestion and the different initial values of $\bi^{dm}$.}\label{AsBint}
	\end{table}
\end{itemize}
\subsubsection*{Direct effect of DF on the function of the small intestine}
Besides the modification of the bolus in the stomach, presence of insoluble DF changes also the physiochemical characteristics of bolus (Section \ref{Model}).

To observe the effect of DF in the model of digestion, value of insoluble and soluble DF was increased from $1\,g$ to $5\, g$ in a bolus of $120\, g$.
 Figure \ref{FsFins} shows that the presence of insoluble DF promotes intestinal absorption, however this increase in absorbed dry nutrients is not meaningful. The results in this figure, show that the increasing of the value of soluble DF decreases the quantity of absorbed dry nutrients and increases the final total mass.
   
 Therefore, numerical simulation shows that the positive effect of insoluble DF on the amount of absorbed dry nutrients is mainly due to the modification of the bolus in the stomach. However, the effect of the interaction of bolus with the small intestine as described in Figure \ref{FsFins} may not be neglected. 
	\begin{figure}[h]
	\subfigure[Absorption versus $F_{ins}$]{
	\includegraphics[width=70mm]{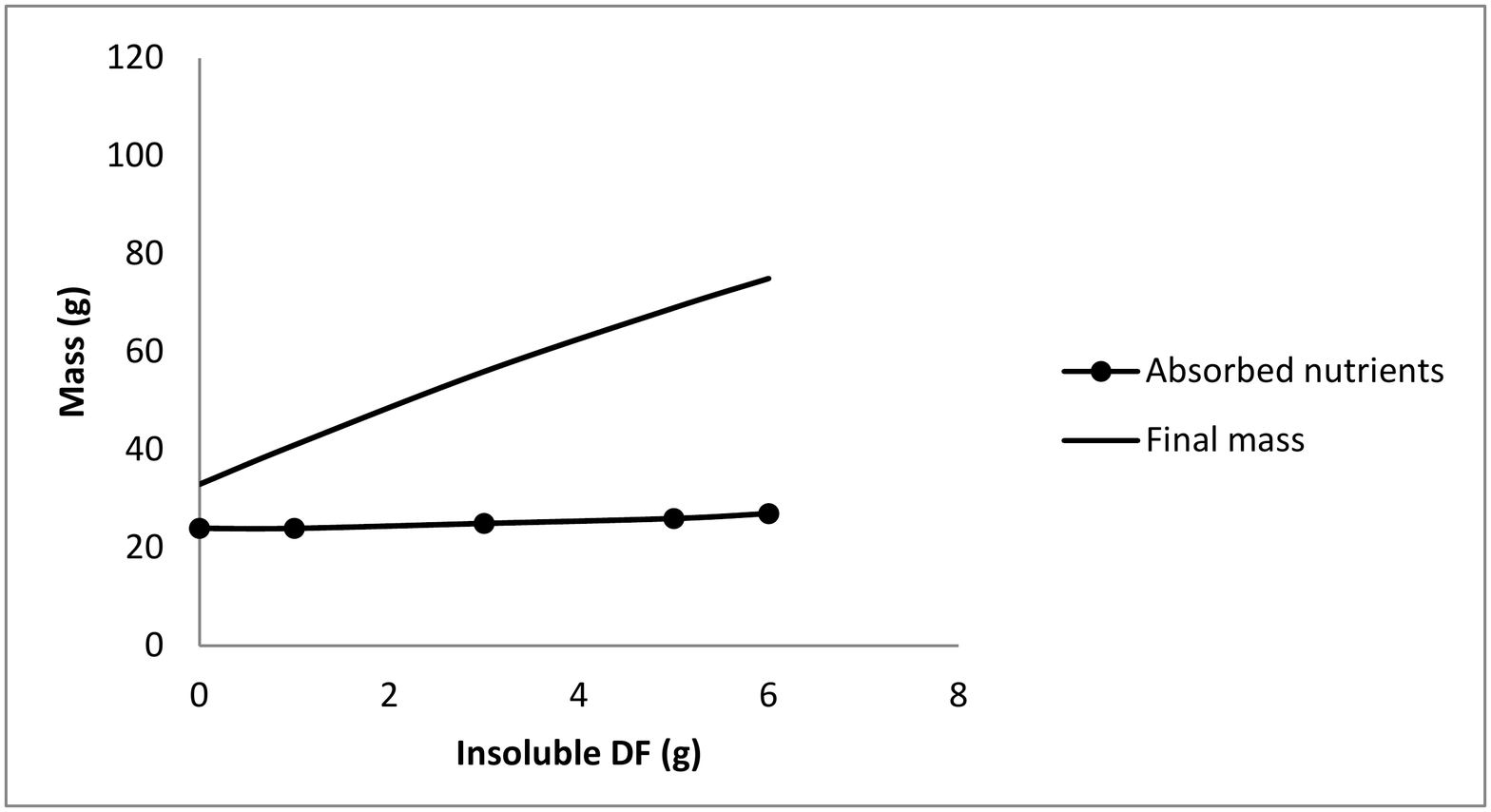} 
	  }
	\subfigure[Absorption versus $F_{s}$]{
	\includegraphics[width=70mm]{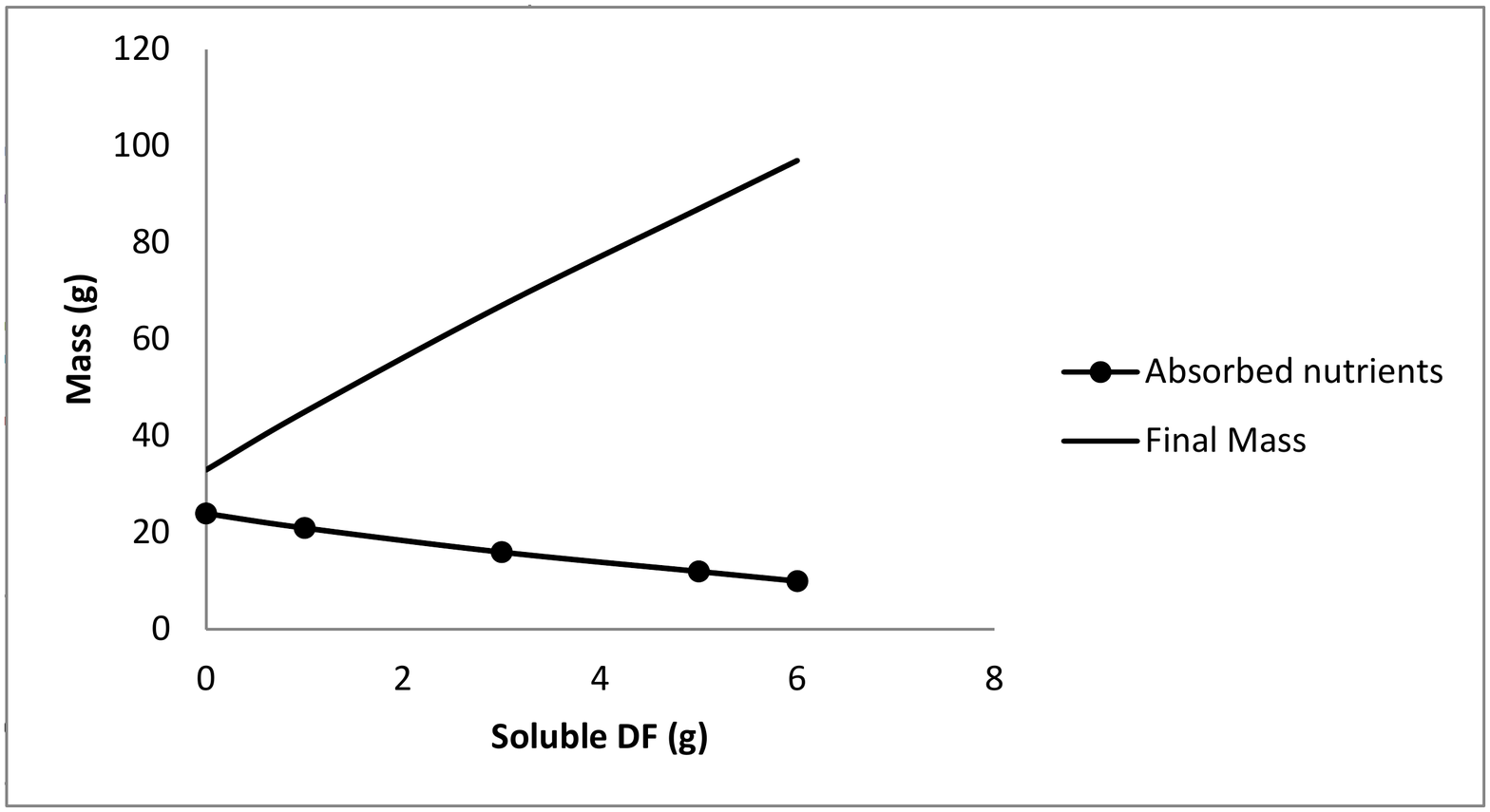} 
	  }
	\caption{\small Change in the final mass of bolus and absorbed dry nutrients for different amount of DF in the bolus at the entry of duodenum.}\label{FsFins}
	\end{figure}
\subsubsection*{Time of intestinal transit in the presence of DF}
Numerical results of transit time in presence of DF are presente in Table \ref{timefsol}. 

These results show that integrating insoluble DF in the bolus decreases the time of intestinal transit from $5\, h$ for a bolus of $120 \, g$ without insoluble DF to $3,9\,h$ for a bolus of the same mass which contains $5\, g$ of DF. These results are consistent with published values. The experiences done by Wilfart et al. \cite{Wilfart2007} have shown that increasing dietary fibre content reduced or tended to reduce the mean retention time in the small intestine.

These results show that integrating soluble DF in the bolus increases the intestinal digestion time from $5\,h$ to $6,7\,h$ illustrating the effect of viscosity due to soluble DF on transit time.

	\begin{table}[h]
	\begin{tabular}{|c|c|c|}
	  \hline
	\multicolumn{1}{|c|}{ $x=0$} &\multicolumn{2}{c|}{Intestinal transit time (hour)}  \\
	\hline
	 $F\cdot DM^{-1}$ (\%) & bolus containing $F_{insol}$  & bolus containing $F_{sol}$ \\
	\hline
	\centering $0$&\centering $5$ & $5$\\
	\centering  $2$&\centering$4,7$ &$5,4$\\  
	 \centering $7$&\centering$4,3$ &$5,9$\\
	\centering  $11$&\centering$4$ &$6,5$\\
	\centering $14$ &\centering$3,9$ &$6,7$ \\
	\hline
	\end{tabular}\caption{\small Intestinal transit time for the different quantities of $F=F_{sol}$ or $F_{insol}$ in the bolus at the entry of duodenum}\label{timefsol}
 	\end{table}
%
\subsection{Water associated to dry matter}\label{testw}
To study the influence of the quantity of associated water on digestion, variation of the values of $\a$, $\b$ and $\g$ have been tested. We observe their influence on digestion and specifically on the absorbed dry nutrients and $\as$-$A_{ns}$ equilibrium. 
Two main simulations are carried out : an uniform water content for $\as$, $\bi$ and $\ba$ i.e. $\a=\b=\g$ and a a non uniform one.

The objective of these experiences is to understand how the different values of $\a,\,\b$ and $\g$ influence the digestion in our modelling.
\subsubsection{ Uniform water content for $\as$, $\bi$ and $\ba$ }
The value of $\a=\b=\g$ varied from $1$ to $4$ in the model presented in section \ref{Model}. Our objective is to observe its effect on the value of absorbed dry nutrients as well as on the final mass of bolus.

Numerical results presented in Figure \ref{water1} show the negative effect of this increase on the absorbed dry nutrients. Increasing the quantity of water ($\a$, $\b$ and $\g$) associated with the dry feedstuffs ($\as^{dm}$, $\bi ^{m}$ and $\ba ^{dm}$) in our model, dilutes the bolus and decreases the volumic degradation, it decreases also the quantity of dry nutrients in contact with the internal surface of the bolus.  
%
	\begin{figure}[h]
	    \centering
	\includegraphics[width=100mm]{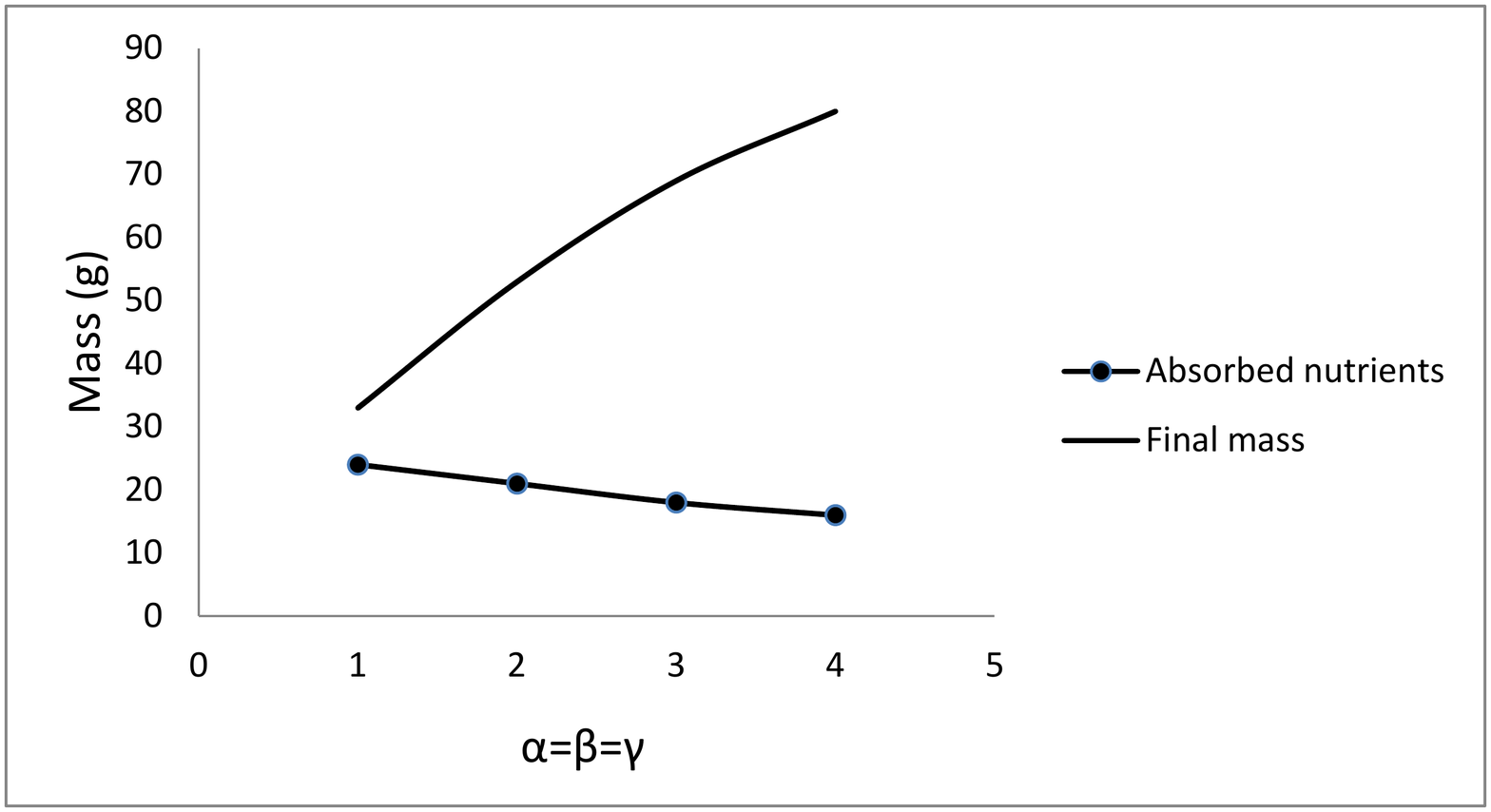}
	\caption{\small Dependence of the absorbed dry nutrients and the final mass of bolus at the end of the small intestine on the value of $\a$, $\b$ and $\g$.}\label{water1}
	\end{figure}
%
These results seems to be consistent with the reality, in fact, the more water is presented in the bolus, the less (pancreatic and brush border) enzymes and molecules are in contact.

We are also interested by the effect of these variations on the $\as^{dm}-A_{ns}$ equilibrium defined by Equation (\ref{equilAns}). The results are shown in Figure \ref{abc1234}.
The equilibrium is almost reached in the four experiences, however the choice of $k_{equi}$
can change the necessary time to reach the equilibrium. 
A significant production of $A_{ns}$ is observed in Figure \ref{abc1234}.$d$ because of the lack of the available water at the beginning of the small intestine.

	\begin{figure}[h]
	\subfigure[$\a=\b=\g=1$]{
	\includegraphics[width=70mm]{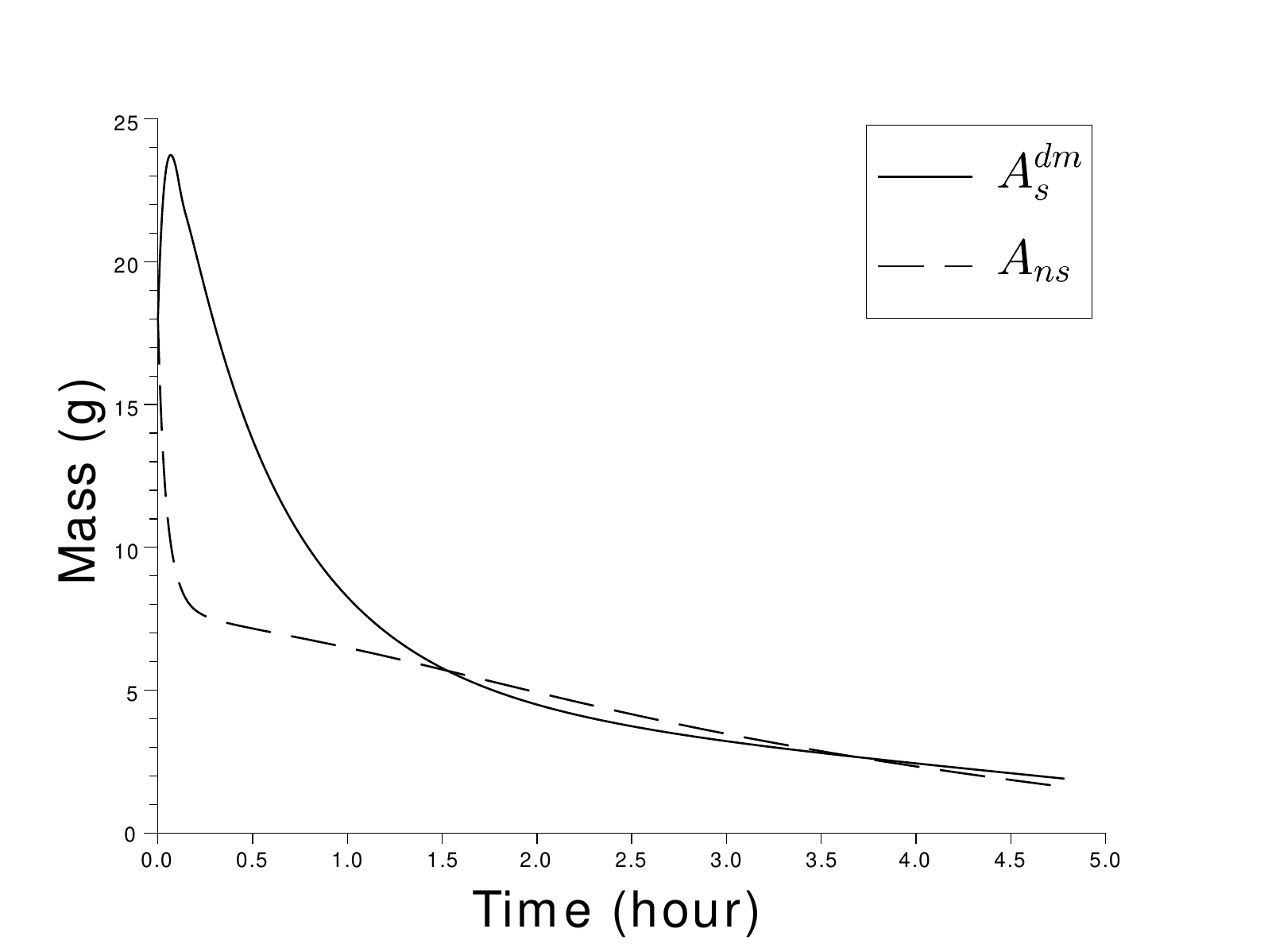} 
	  }
	\subfigure[$\a=\b=\g=2$]{
	\includegraphics[width=70mm]{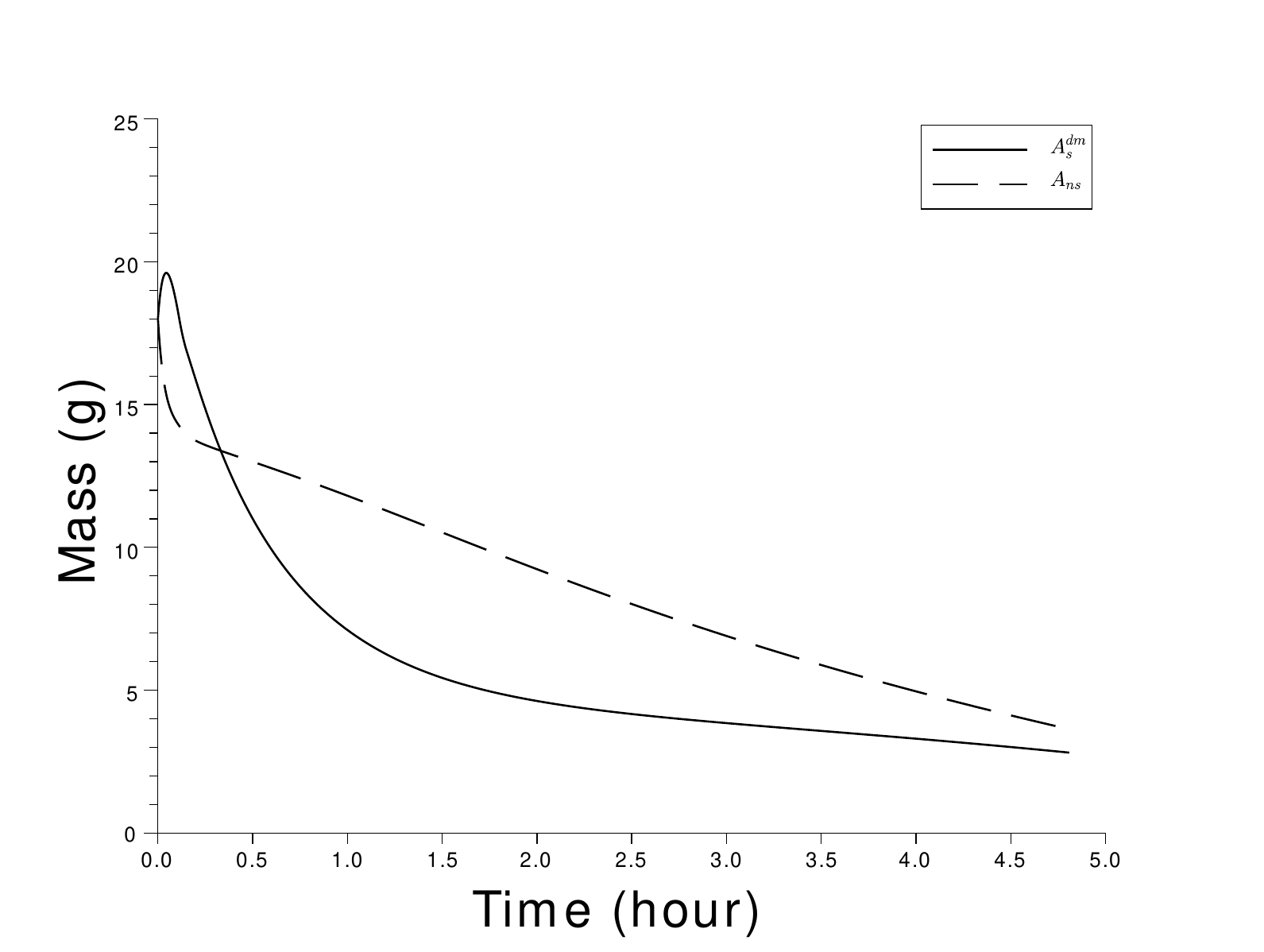} 
	  }
	\subfigure[$\a=\b=\g=3$]{
	\includegraphics[width=70mm]{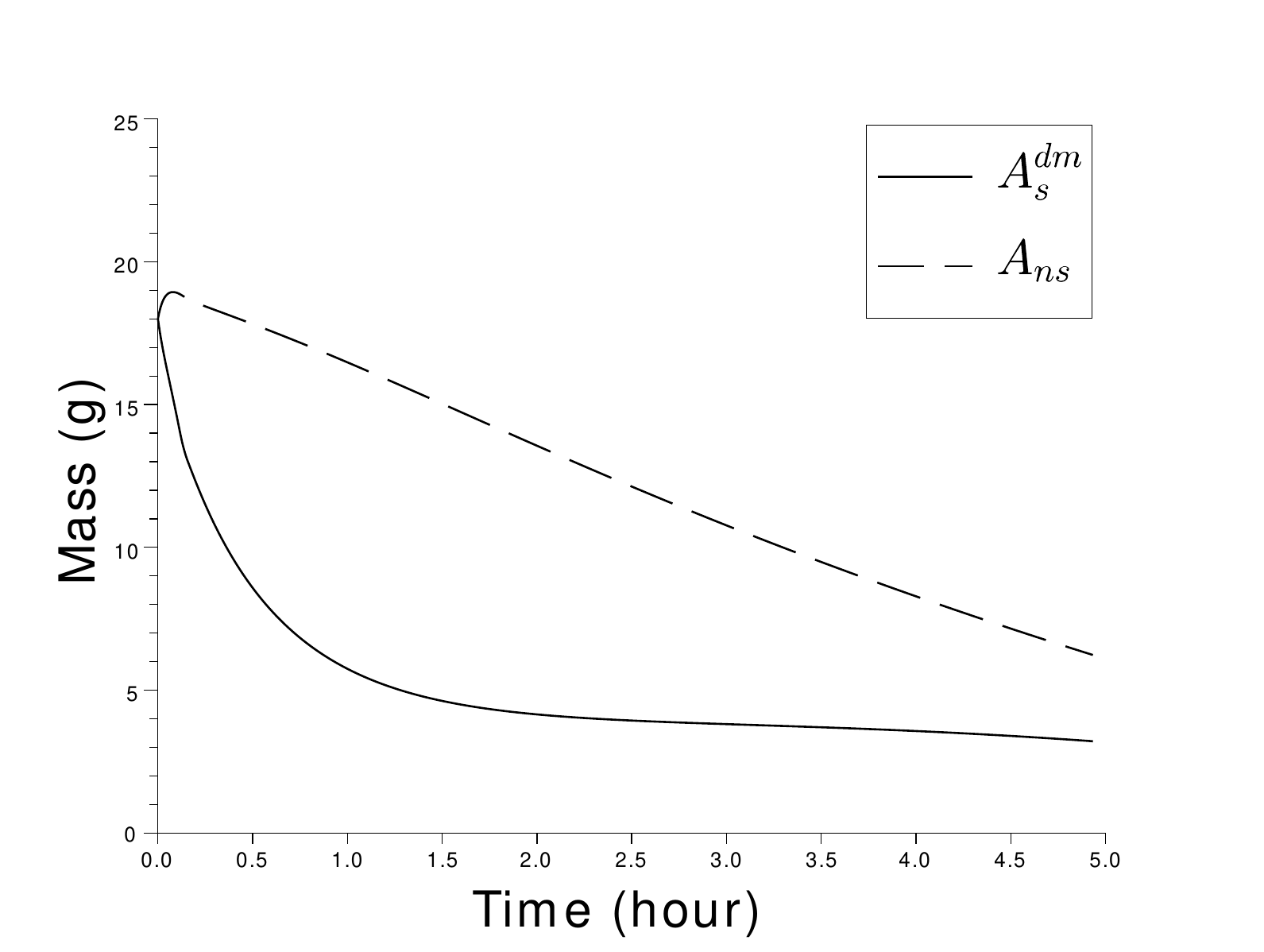}
	  }
	\subfigure[$\a=\b=\g=4$]{
	\includegraphics[width=70mm]{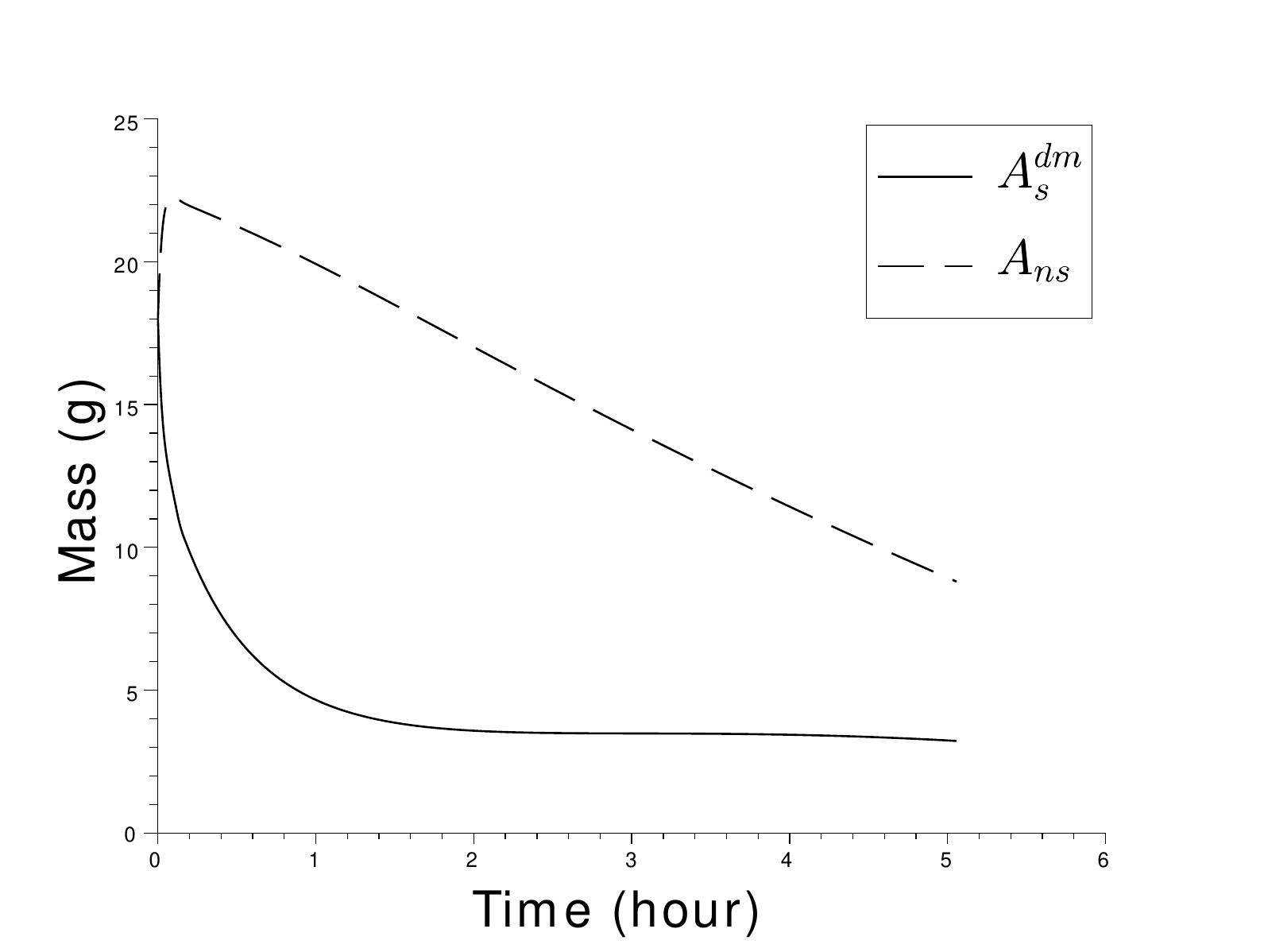}
	  }
	\label{equilfig}
	\caption{\small The evolution of the equilibrium $\as^{dm}-A_{ns}$ depends strongly on the quantity of $\a$, $\b$ and $\g$. }\label{abc1234}
	\end{figure}
\subsubsection{Non-uniform water content for $\as$, $\bi$ and $\ba$ }
In the second experience, the quantity of absorbed dry nutrients at the end of the small intestine and the numerical results of $\as$-$A_{ns}$ equilibrium for different values of $\a$, $\b$ and $\g$ were observed. The choice of the values of $\as$, $\bi$ and $\ba$ is based on the hypothesis that the value of $\b$ is always between the values of $\a$ and $\g$ because of the molecule size of $\bi^{dm}$.
	\begin{table}[h]
	\centering%
	\begin{tabularx}{\linewidth}{|X|X|X|}
	  \hline
	&\small Absorbed dry nutrients to $DM$ ratio (\%)&\small Retention time in the small intestine (h)  \\
	\hline
	\small $\a=1,\, \b=2,\,\g=3$ & $45$ & $5,2$
 \\
	\hline
	\small $\a=2\,\b=2,\,\g=2$ &  $50$ &$4,9 $  \\
	\hline
	\small $\a=3,\,\b=2,\,\g=1$& $54$ &$4,6$\\
	\hline
	\end{tabularx}\caption{\small Transit time and absorbed dry nutrients at the end of the small intestine depend on the values of $\a$, $\b$ and $\g$.}\label{abctab}
	\end{table}
Even if the longest transit time was observed for $\a<\b<\g$, it shows the lowest level of dry  absorbed nutrients.
This stresses the key-role of ``available water'' $W$ on digestion. 

$\as-A_{ns}$ equilibrium is observed in two cases : 
\begin{itemize}
 \item[$(i)$] 
Firstly, it was assumed that $\a>\b>\g$. There is a little production of $A_{ns}$ from $\as$ (desoubilization) because of the lack of water at the entry of the small intestine followed by the degradation of $\as$.  

 \item[$(ii)$] 
Secondly,  it was assumed that $\a<\b<\g$. There is a rapid solubilization of $A_{ns}$ followed by degradation of $\as$ into $\bi$. The released water by degradation of $\as$ is not sufficient to maintain the solution state of $\bi$ since $\a<\b$. The required water is provided by the transformation of $\as$ into $A_{ns}$ as it can be seen in the Figure \ref{water3}.
	\begin{figure}[h]
	\subfigure[$\a>\b>\g$]{
	\includegraphics[width=70mm]{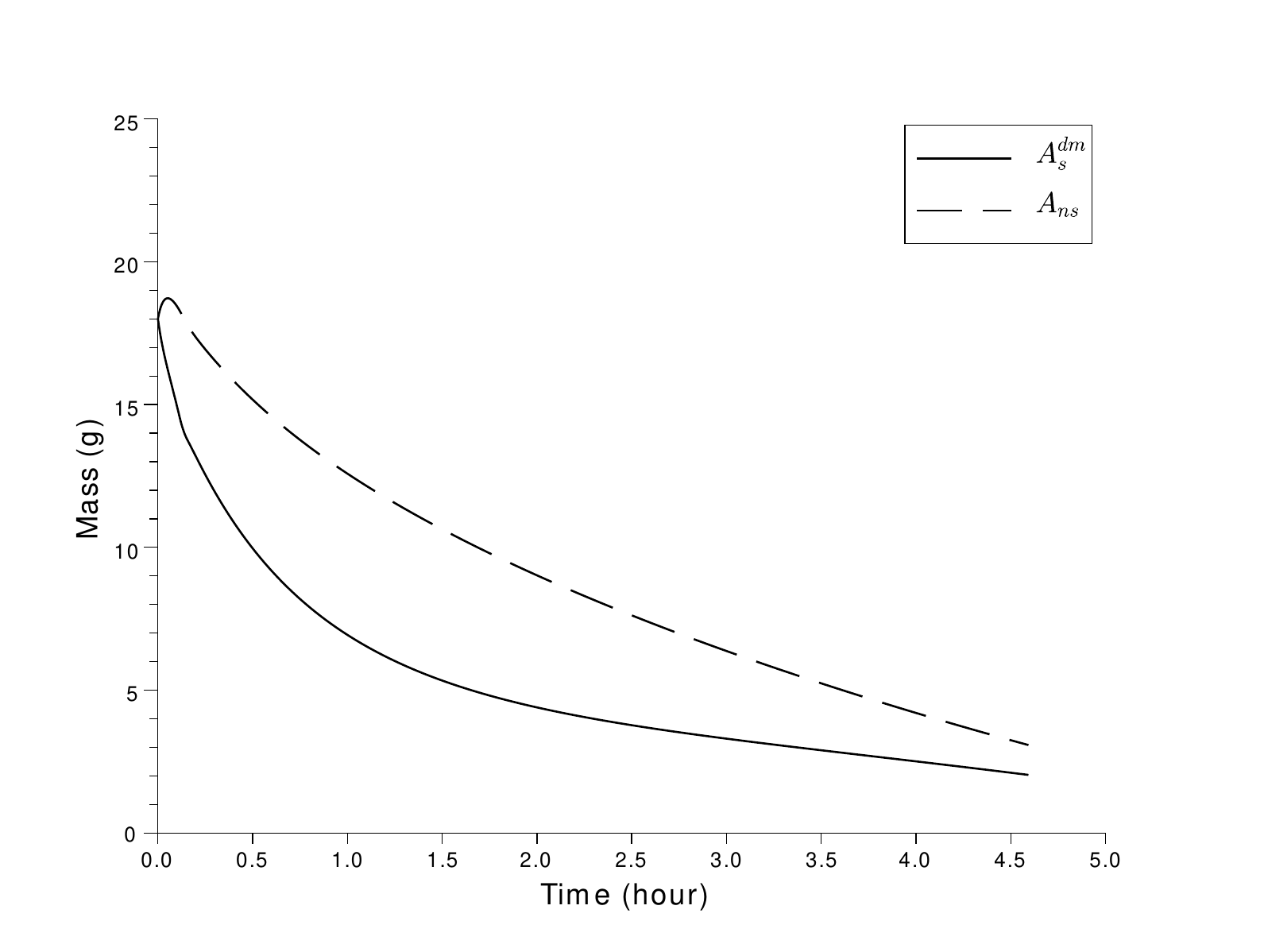} 
	  }
	\subfigure[$\a<\b<\g$]{
	\includegraphics[width=70mm]{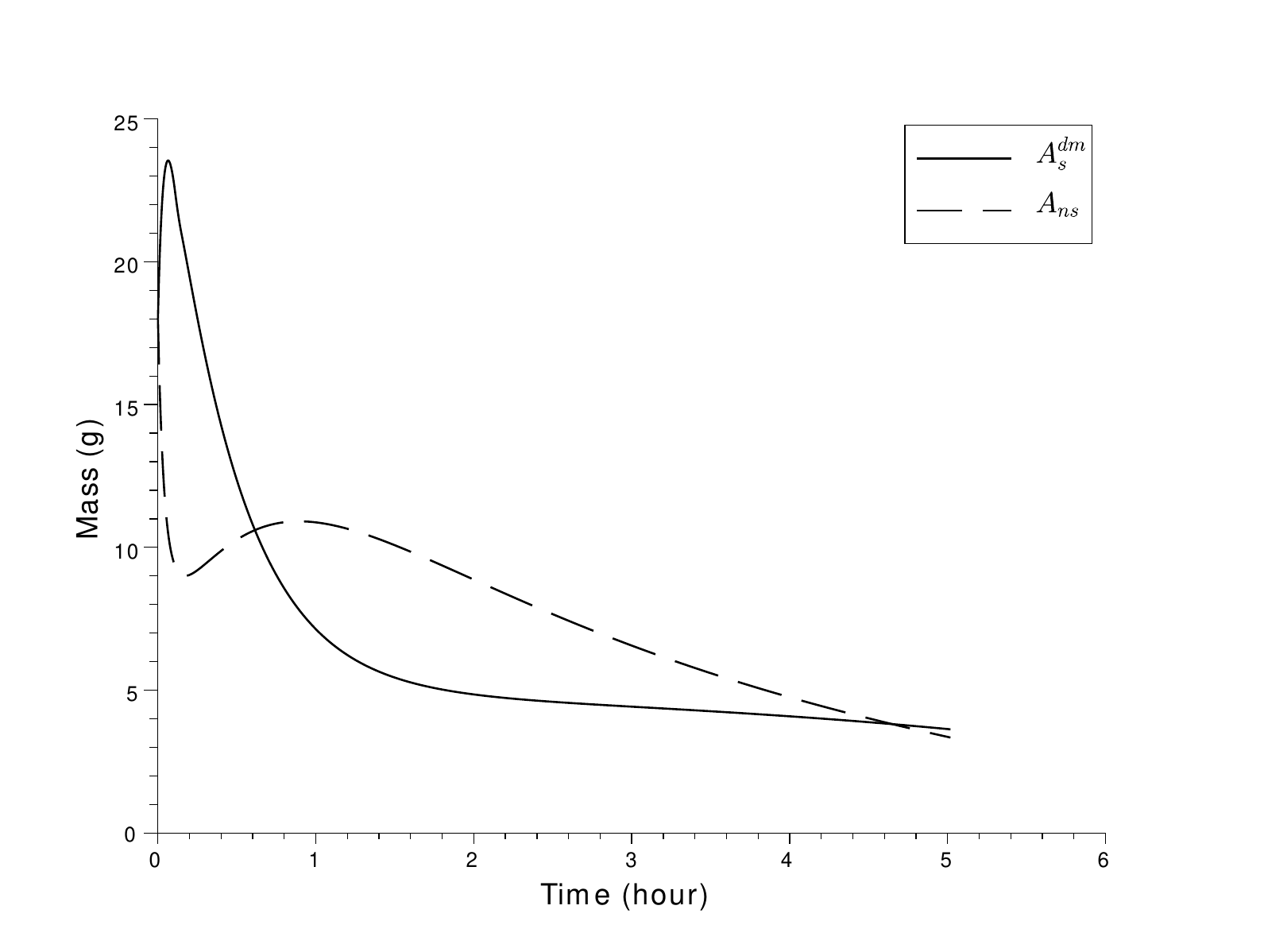}
	  }
	\label{equilfig}
	\caption{\small The equilibrium $\as^{dm}-A_{ns}$ for different quantity of $\a$, $\b$ and $\g$. } \label{water3}
	\end{figure}
 \end{itemize}
\subsection{The variation of the water ratio inside the bolus}\label{testwv}
We are interested by the change in the ratio of $W$ inside the bolus and its influence on the absorption in the small intestine. For the sake of simplicity, it was assumed that $W_{feed}=W_{sec}=0$. 

In order to study this effect, the value of dry matter in the bolus is assumed fixed at $42\, g$ while the value of water increases in each experience. 

In the numerical simulations, the ratio of water included in the bolus represented $55\%$, $60\%$ and $66\%$ of bolus. The \textit{ratio of absorbed dry nutrients to the total absorption (water and nutrients)} is collected. Results are presented in Table \ref{WvsDM}.
	\begin{table}[h]
	\centering%
	\begin{tabularx}{\linewidth}{|X|X|X|c|}
	  \hline
	 \centering $(\%)$ Amount of water at $x~=~0$ inside bolus  &\centering $(\%)$ Ratio of absorbed dry nutrients at $x=17$ &\centering $(g)$ Mass of Absorbed dry nutrients&Transit time (hour)\\
	\hline
	 \centering$50$ &\centering$50$&\centering$25$&$5,7$ \\ 
	\centering$60$ & \centering$35$ &\centering$26$&$5,2$\\
	\centering$66$ & \centering$28$ &\centering$27$&$4,5$ \\
	\hline
	\end{tabularx}\caption{\small Dependence of absorbed dry nutrients to the water ratio. DM stands for the amount of dry matter in the bolus at the entry of bolus. }\label{WvsDM}
	\end{table}
The numerical results of Table \ref{WvsDM} shows that increasing the value of water in the bolus decreases the \textit{ratio of absorbed dry nutrient to the total absorbed matter (water+dry nutrients)}, even though the value of absorbed dry nutrients increases.

 In fact, in our model, at the end of each experience the value of $W_t$ is approximatively $55\%$ of the total mass of bolus and this equilibrium is achieved almost quickly (because of the choice of $k_w$). Therefore, increasing the value of water in this model does not have a meaningful effect on the final absorption and the slight increase in the absorbed mass of dry nutrients is due to the change of the volume of bolus in each experience which promotes the access of nutrients to intestinal wall for absorption, although this increase dilutes the bolus and decreases the volumic degradation. Here again, this is the direct result the choice of the parameters (rate of surfacic $(k_{surf})$ and volumic $(k_{vol})$ degradation ). 
\section{Discussion}\label{disc}
In this paper, we have continued the modelling of the  digestion in the small intestine started in \cite{Taghipoor2011}. The objective was to obtain a more realistic model of digestion process by including new phenomena and completing the others:   
$DM$ and water are treated separately, water evolution is influenced by the degradation of other molecules, the effects of DF on the digestion are taken into account which is also a first step toward a non-homogeneous model with different types of feedstuffs. 

The advantages and the defects of this model as  well as the perspectives are outlined in the following paragraphs.

One of the main aspect of this model remains its genericity, we have tried to identify and model the main generic phenomena of digestion and ignored or implicitly taken into account in the parameters the ones which required the specific properties of feedstuffs' molecules (effects on gastric emptying, on viscosity,...etc.).   
The different steps of digestion (equilibrium between $A_{ns}$ and $ A_s$, successive transformation of $\as$ into $B_{int}$ and $\ba $), the effects of physical characteristic of bolus (surface and volume) on its degradation and transit, the interaction between DF and feedstuffs molecules (among others)  have been considered while some other phenomena like the separation between the enzyme activities of different feedstuffs' molecules, dependence of the enzyme activity with respect to the dilution of bolus, the different substrates density, the impact of DF on initial condition (and others) have been ignored. 
To our point of view, we obtained a more realistic model by integrating these new phenomena in the model of digestion. In particular, the WHC of dietary fibre which in turn interferes the digestion of other feedstuffs molecules, leads us to introduce the separation between  DM and water and consider all the effects of water.

However, some other phenomena are still ignored either because of the lack of information concerning their effects or because of their supposed little impact on digestion at this scale. 
Of course, It would be interesting to include in the model the phenomena like interaction between different categories of feedstuffs and then define the specific enzyme activity for different cases to exploit their potential to impact the digestion. The future development of the model will be based on these new objectives. 
Modelling the influence of soluble and insoluble DF on the initial condition and on degradation of other feedstuffs' molecules as well as on the movement of bolus (Experience \ref{testf}) is the first try to model a more realistic non-homogeneous bolus.       
DF have normally a high WHC which increases the volume of aqueous phase in the bolus and therefore dilutes the solution of nutrients and enzymes \cite{O.1998}. This is known to influence the volumic and surfacic reactions. However, these effects depend highly to the choice of $W_{sol}$ and $W_{insol}$. Another aspect of the model is the water equilibrium and its impact on absorption, here again the choice of parameters ($W_{as}$, $W_{int}$ and $W_{abs}$) plays significant role on the final results of digestion.
On the other hand, Experience \ref{testwv} reveals the role of other parameters of model ($k_w$, $k_{equi}$, $k_{surf}$, $k_{vol}$). As it has been described in this experience, the change of water absorption rate $k_w$, can change ( even inverse) the numerical results. However, it is worth pointing that the choice of the model parameters are based on the observed behaviours (literature), the results are therefore consistent qualitatively with the reality (positive influence of insoluble DF on digestion, negative effect of soluble DF, ...). 
Taking into account these new phenomena requires the introduction of new parameters which can be identified with the help of existing experimental data. To the best of our knowledge, some of these parameters like the water associated to dry nutrients ($W_{sol}$ , $W_{insol}$, $W_{as}$, $W_{int}$ and $W_{abs}$) are introduced for the first time and they should be identified

One of the advantages of this model is its capacity to be reduced and to be adapt to the new experiences which makes the parameters identification possible. Reduction consists in the decreasing the number of equations of system or the number of parameters without affecting its genericity (e.g. a bolus which does not contain the DF, results in a more simplified digestion process which in turn caused decreasing the equations (parameters) involved in the digestion model).

It is also worth pointing that the value of most of the parameters depends to the special category of feedstuffs. The close collaborations between biologists and mathematicians is therefore required to identify these new parameters (literature data in biology, define the new experimentations, ...). This reveals one of the main interest of modelling which is to ask the precise questions about the modeled phenomenon. In fact, this approach allows to use all the existing data and limit the new animal experimentation to the special cases (when the existing data are not sufficient).    

 \section{Appendix}
\subsection{The Model equations}
\begin{itemize}
 \item[$\diamond$] Equation of transport
		\begin{align*}\label{TransISolnSolDF}
		\dfrac{d^2x}{dt^2}=\tau(1-C^{-1}\frac{dx}{dt})\frac{c_0+c_1r_{sol}}{a+bx}-\dfrac{K_{visco}}{[W]}\frac{dx}{dt}.                 
		\end{align*}
\item[$\diamond$] Non solubilized substrate  $A_{ns}$
		\begin{equation*}
		\dfrac{d{A_{ns}}}{dt}=-k_{equi}\big(\mu([W])A_{ns}-\as\big)
		\end{equation*}
 \item[$\diamond$] Soluble DF $F_{sol}$ 
		\begin{equation*}\label{Fs}
		\dfrac{d{F}_{sol}^{dm}}{dt}= -k_{s}e_{exo}\tilde{ph}(x)[F_{sol}^{dm}]V_{app}
		\end{equation*}
\item[$\diamond$] Solubilized substrate $A^{dm}_s$
		\begin{align*}
		\dfrac{d{A_s^{dm}}}{dt}=-k_{vol}(x)[\as^{dm}]V_{app}-k_{surf}[\as][W]S_{sol}
		+k_{equi}\big(\mu([W])A_{ns}-\as\big)
		+\hbox{secretions}
		\end{align*}
\item[$\diamond$] Intermediate substrate $B^{dm}_{int}$
		\begin{align*}
		\dfrac{d{\bi^{dm}}}{dt}=&k_{vol}(t)[\as^{dm}]V_{app}+k_{s}e_{exo}\tilde{ph}(x)[F_{sol}^{dm}]V_{app}- \tilde{k}_{surf}[\bi][W]S_{sol}+\hbox{secretions}
		\end{align*}
\item[$\diamond$] Absorbable nutrients $\ba^{dm}$
		\begin{align*}
		\dfrac{d{\ba^{dm}}}{dt}=\bigg(\tilde{k}_{surf}[\bi]+k_{surf}[\as]  \bigg) [W]S_{sol}- k_{abs}[\ba^{dm}]S_{sol}
		\end{align*}
 \item[$\diamond$] Exogenous enzymes $e(t)$ 
		\begin{equation*}\label{Fs}
		\dfrac{d{e}}{dt}= -k_{e}e
		\end{equation*}
 \item[$\diamond$] Free Water Evolution  W(t)
		\begin{align*}\label{equW}
		\dfrac{d{W}}{dt}=&-\a k_{equi}\big(\mu([W])A_{ns}-\as\big)+(\a-\b) k_{vol}(t)[\as^{dm}]V_{app}\nonumber\\
		&+[W]\bigg( (\b-\g)\tilde{k}_{surf}[\bi]+ (\a-\g)k_{surf}[\as]\bigg)\nonumber\\
		&+(\lambda_s-\beta)k_{s}e_{exo}\tilde{ph}(x)[F_{sol}^{dm}]V_{app}+\g k_{abs}[\ba^{dm}]S_{sol}-k_w(W-0.1M).
		\end{align*}
 \item[$\diamond$] Evolution of the bolus mass M(t) 
		\begin{equation*}
		\dfrac{d{M}}{dt}= \dfrac{M}{M-W}\big(-k_w(W-0.1M)- k_{abs}[\ba^{dm}]S_{sol} + \hbox{Secretions}\big)
		\end{equation*}
 \item[$\diamond$] Evolution of bolus volume V(t) 
	    \begin{equation*}
	      \dfrac{d{V}}{dt}= -k_w(W-0.1M)+ \hbox{Secretions}
	    \end{equation*}

\end{itemize}

%

\bibliographystyle{unsrt}
\bibliography{bibliog}

\end{document}